\documentclass{article}
\usepackage{amsmath,amsthm,amssymb,amsfonts}
\usepackage{times}
\usepackage{enumerate}
\usepackage{amssymb}
\usepackage{graphicx}
\usepackage[noadjust]{cite}
\usepackage{color}

\newtheorem{thm}{Theorem}[section]
\newtheorem{cor}[thm]{Corollary}
\newtheorem{lem}[thm]{Lemma}

\newtheorem{prop}[thm]{Proposition}

\newtheorem{rem}[thm]{Remark}
\newtheorem{example}[thm]{Example}
\makeatletter
\def\author#1{\gdef\autrun{\def\and{\unskip, }#1}\gdef\@author{#1}}

\makeatother

\makeatletter
\let\@fnsymbol\@alph
\makeatother

\DeclareMathOperator*{\esssup}{ess\,sup}
\DeclareMathOperator*{\essinf}{ess\,inf}
\DeclareMathOperator*{\linspan}{span}

\begin{document}

\title{Compactness and existence results for the $p$-Laplace equation}

\author{Marino Badiale\thanks{Dipartimento di Matematica ``Giuseppe Peano'', Universit\`{a} degli Studi di
Torino, Via Carlo Alberto 10, 10123 Torino, Italy. 
e-mails: \texttt{marino.badiale@unito.it}, \texttt{michela.guida@unito.it}}
\textsuperscript{,}\thanks{Partially supported by the PRIN2012 grant ``Aspetti variazionali e
perturbativi nei problemi di.renziali nonlineari''.}
\ -\ Michela Guida\textsuperscript{a,}\thanks{Member of the Gruppo Nazionale di Alta Matematica (INdAM).}
\ -\ Sergio Rolando\thanks{Dipartimento di Matematica e Applicazioni, Universit\`{a} di Milano-Bicocca,
Via Roberto Cozzi 53, 20125 Milano, Italy. e-mail: \texttt{sergio.rolando@unito.it}}%
\ \textsuperscript{,\,c}
}

\date{}
\maketitle

\begin{abstract}
Given $1<p<N$ and two measurable functions $V\left( r\right) \geq 0$ and $%
K\left( r\right) >0$, $r>0$, we define the weighted spaces 
\[
W=\left\{ u\in D^{1,p}(\mathbb{R}^{N}):\int_{\mathbb{R}^{N}}V\left( \left|
x\right| \right) \left| u\right| ^{p}dx<\infty \right\} ,\quad
L_{K}^{q}=L^{q}(\mathbb{R}^{N},K\left( \left| x\right| \right) dx) 
\]
and study the compact embeddings of the radial subspace of $W$ into $%
L_{K}^{q_{1}}+L_{K}^{q_{2}}$, and thus into $L_{K}^{q}$ ($%
=L_{K}^{q}+L_{K}^{q}$) as a particular case. 
We consider exponents $q_{1},q_{2},q$ that can be greater or smaller than $p$.
Our results do not require any compatibility between how the potentials $V$
and $K$ behave at the origin and at infinity, and essentially rely on power
type estimates of their relative growth, not of the potentials separately.
We then apply these results to the investigation of existence and
multiplicity of finite energy solutions to nonlinear $p$-Laplace equations
of the form 
\[
-\triangle _{p}u+V\left( \left| x\right| \right) |u|^{p-1}u=g\left( \left|
x\right| ,u\right) \quad \text{in }\mathbb{R}^{N},\ 1<p<N, 
\]
where $V$ and $g\left( \left| \cdot \right| ,u\right) $ with $u$ fixed may
be vanishing or unbounded at zero or at infinity. Both the cases of $g$
super and sub $p$-linear in $u$ are studied and, in the sub $p$-linear case,
nonlinearities with $g\left( \left| \cdot \right| ,0\right) \neq 0$ are also
considered.\bigskip

 \smallskip

 \medskip
\end{abstract}


\section{Introduction}

In this paper we pursue the work we made in papers \cite{BGR_I,BGR_II,GR-nls}%
, where we studied embedding and compactness results for weighted Sobolev
spaces in order to get existence and multiplicity results for semilinear
elliptic equations in $\mathbb{R}^{N}$, by variational methods. 

In the present paper, we face nonlinear elliptic $p$-Laplace equations with
radial potentials, whose prototype is 
\begin{equation}
-\triangle _{p}u+V\left( \left| x\right| \right) \left| u\right|
^{p-1}u=K\left( \left| x\right| \right) f\left( u\right) \quad \text{in }%
\mathbb{R}^{N}  \label{EQ}
\end{equation}
(more general nonlinear terms will be actually considered in the following).
Here $1<p<N$, $f:\mathbb{R}\rightarrow \mathbb{R}$ is a continuous
nonlinearity satisfying $f\left( 0\right) =0$ and $V\geq 0,K>0$ are given
potentials, which may be vanishing or unbounded at the origin or at infinity.

To study this problem we introduce the weighted Sobolev space 
\[
W:=\left\{ u\in D^{1,p}(\mathbb{R}^{N}):\int_{\mathbb{R}^{N}}V\left( \left|
x\right| \right) \left| u\right| ^{p}dx<\infty \right\} 
\]
equipped with the standard norm 
\begin{equation}
\left\| u\right\| ^{p}:=\int_{\mathbb{R}^{N}}\left( \left| \nabla u\right|
^{p}+V\left( \left| x\right| \right) \left| u\right| ^{p}\right) dx,
\label{norm}
\end{equation}
%
%
%
and we say that $u\in W$ is a \textit{weak solution}\emph{\ }to (\ref{EQ})
if 
\[
\int_{\mathbb{R}^{N}}|\nabla u|^{p-2}\nabla u\cdot \nabla h\,dx+\int_{%
\mathbb{R}^{N}}V\left( \left| x\right| \right) \left| u\right|
^{p-2}uh\,dx=\int_{\mathbb{R}^{N}}K\left( \left| x\right| \right) f\left(
u\right) h\,dx\quad \text{for all }h\in W. 
\]
The natural approach in studying weak solutions to equation (\ref{EQ}) is
variational, since these solutions are (at least formally) critical points
of the Euler functional 
\[
J\left( u\right) =\frac{1}{p}\left\| u\right\| ^{p}-\int_{\mathbb{R}%
^{N}}K\left( \left| x\right| \right) F\left( u\right) dx, 
\]
where $F\left( t\right) :=\int_{0}^{t}f\left( s\right) ds$. Then the problem
of existence is easily solved if $V$ does not vanish at infinity and $K$ is
bounded, because standard embeddings theorems of $W$ and its radial subspace 
$W_{r}$ into the weighted Lebesgue space 
\[
L_{K}^{q}:=L_{K}^{q}(\mathbb{R}^{N}):=L^{q}(\mathbb{R}^{N},K\left( \left|
x\right| \right) dx) 
\]
are available (for suitable $q$'s). As we let $V$ and $K$ to vanish, or to
go to infinity, as $|x|\rightarrow 0$ or $|x|\rightarrow +\infty $, the
usual embeddings theorems for Sobolev spaces are not available anymore, and
new embedding theorems need to be proved. This has been done in several
papers: see e.g. the references in \cite{BGR_I,BGR_II,GR-nls} for a
bibliography concerning the usual Laplace equation, and \cite
{Anoop,Su12,Cai-Su-Sun,SuTian12,Su-Wang-Will-p,Yang-Zhang,Zhang13,BPR,Zhao-Su-Wang-16,Bart-Candela-Salv-16,Chen-Wang-15,Li-Cai-Su-14}
for equations involving the $p$-laplacian.

The main novelty of our approach (in \cite{BGR_I,BGR_II} and in the present
paper) is two-fold. First, we look for embeddings of $W_{r}$ not into a
single Lebesgue space $L_{K}^{q}$ 
but into a sum of Lebesgue spaces $L_{K}^{q_{1}}+L_{K}^{q_{2}}$. This allows
to study separately the behaviour of the potentials $V$ and $K$ at $0$ and $%
\infty $, and to assume different set of hypotheses about these behaviours.
Second, we assume hypotheses not on $V$ and $K$ separately but on their
ratio, so allowing asymptotic behaviours of general kind for the two
potentials.

Thanks to these novelties, our embedding results yield existence of
solutions for (\ref{EQ}) in cases which are not covered by the previous
literature. Moreover, one can check that our embeddings are also new in some
of the cases already treated in previous papers (see e.g. Example \ref
{EX(ST)}), thus giving existence results which improve some well-known
theorems in the literature. %
%

The proofs of our embedding theorems for the space $W_{r}$ are
generalizations of those presented in \cite{BGR_I} for the Hilbertian case $%
p=2$. The generalizations to the case $1<p<N$ are not difficult but boring
and lengthy, because one needs to repeat a lot of detailed computations, the
basic ideas remaining the same. In view of this, in the present paper we
limit ourselves to state our embedding results and to present in detail some
examples, leading to new existence results for equation (\ref{EQ}). For all
the proofs, with full details, we refer the reader to the specific document 
\cite{BGR-p.arxiv}, which is essentially a longer version of Section \ref
{SEC: comp} below.

This paper is organized as follows. In Section \ref{SEC: comp} we state our
main results: a general result concerning the embedding properties of $W_{r}$
into $L_{K}^{q_{1}}+L_{K}^{q_{2}}$ (Theorem \ref{THM(cpt)}) and some
explicit conditions ensuring that the embedding is compact (Theorems \ref
{THM0}, \ref{THM1}, \ref{THM2} and \ref{THM3}). In Section \ref{SEC:
examples} we apply our compactness results to some examples, with a view to
both illustrate how to use them in concrete cases and to compare them with
the related literature. In Section \ref{SEC: ex} we present existence and
multiplicity results for equations like (\ref{EQ}), but with more general
nonlinearities, whose proofs are given in Section \ref{SEC:pf-general}%
.\medskip

\noindent \textbf{Notations. } We end this introductory section by
collecting some notations used in the paper. \smallskip

\noindent $\bullet $ We denote $\mathbb{R}_{+}:=\left( 0,+\infty \right) $, $%
B_{R}:=\left\{ x\in \mathbb{R}^{N}:\left| x\right| <R\right\} $, $R>0$, and $%
A^{c}:=\mathbb{R}^{N}\setminus A$ for any $A\subseteq \mathbb{R}^{N}$.


\noindent $\bullet $ $\left\| \cdot \right\| _{X}$ and $X^{\prime }$ denote
the norm and the dual space of a Banach space $X$, in which $\rightarrow $
and $\rightharpoonup $ mean \emph{strong} and \emph{weak }convergence
respectively.

\noindent $\bullet $ $\hookrightarrow $ denotes \emph{continuous} embeddings.

\noindent $\bullet $ $C_{\mathrm{c}}^{\infty }(\Omega )$ is the space of the
infinitely differentiable real functions with compact support in $\Omega
\subseteq \mathbb{R}^{d}$ open. 

\noindent $\bullet $ For any measurable set $A\subseteq \mathbb{R}^{d}$, $%
L^{q}(A)$ and $L_{\mathrm{loc}}^{q}(A)$ are the usual real Lebesgue spaces
and, if $\rho :A\rightarrow \mathbb{R}_{+}$ is a measurable function, $%
L^{p}(A,\rho \left( z\right) dz)$ is the real Lebesgue space with respect to
the measure $\rho \left( z\right) dz$ ($dz$ stands for the Lebesgue measure
on $\mathbb{R}^{d}$). In particular, if $K:\mathbb{R}_{+}\rightarrow \mathbb{R}%
_{+} $ is measurable, we denote $L_{K}^{q}\left( E\right) :=L^{q}\left(
E,K\left( \left| x\right| \right) dx\right) $ for any measurable set $%
E\subseteq \mathbb{R}^{N}$.

\noindent $\bullet $ For $1<p<N$, $p^{*}:=pN/\left( N-p\right) $ is the
Sobolev critical exponent and $D^{1,p}(\mathbb{R}^{N})=\{u\in L^{p^{*}}(%
\mathbb{R}^{N}):\nabla u\in L^{p}(\mathbb{R}^{N})\}$ 
is the usual Sobolev space, which identifies with the completion of $C_{%
\mathrm{c}}^{\infty }(\mathbb{R}^{N})$ with respect to the norm of the
gradient. $D_{\mathrm{rad}}^{1,p}(\mathbb{R}^{N})$ is the radial subspace of 
$D^{1,p}(\mathbb{R}^{N})$.


\section{Compactness results \label{SEC: comp}}

\noindent Assume $1<p<N$ and consider two functions $V,K$ such that:

\begin{itemize}
\item[$\left( \mathbf{V}\right) $]  $V:\mathbb{R}_{+}\rightarrow \left[
0,+\infty \right] $ is a measurable function such that $V\in L^{1}\left(
\left( r_{1},r_{2}\right) \right) $ for some $r_{2}>r_{1}>0;$

\item[$\left( \mathbf{K}\right) $]  $K:\mathbb{R}_{+}\rightarrow \mathbb{R}_{+}$
is a measurable function such that $K\in L_{\mathrm{loc}}^{s}\left( \mathbb{R}%
_{+}\right) $ for some $s>1$. 
\end{itemize}

\noindent Define the function spaces 
\begin{equation}
W:=D^{1,p}(\mathbb{R}^{N})\cap L^{p}(\mathbb{R}^{N},V(|x|)dx),\quad
W_{r}:=D_{\mathrm{rad}}^{1,p}(\mathbb{R}^{N})\cap L^{p}(\mathbb{R}%
^{N},V(|x|)dx)  \label{spaces}
\end{equation}
and let $||u||$ be the standard norm (\ref{norm}) in $W$ (and $W_{r}$).
Assumption $\left( \mathbf{V}\right) $ implies that the spaces $W$ and $%
W_{r} $ are nontrivial, while hypothesis $\left( \mathbf{K}\right) $ ensures
that $W_{r}$ is compactly embedded into the weighted Lebesgue space $%
L_{K}^{q}(B_{R}\setminus B_{r})$ for every $1<q<\infty $ and $R>r>0$ (see 
\cite[Lemma 3.1]{BGR-p.arxiv}). In what follows, the summability assumptions
in $\left( \mathbf{V}\right) $ and $\left( \mathbf{K}\right) $ will not play
any other role than this.

Given $V$ and $K$, we define the following functions of $R>0$ and $q>1$: 
\begin{eqnarray}
\mathcal{S}_{0}\left( q,R\right)&:=& \sup_{u\in W_r,\,\left\| u\right\| =1
}\int_{B_{R}}K\left( \left| x\right| \right) \left| u\right| ^{q}dx,
\label{S_o :=} \\
\mathcal{S}_{\infty }\left( q,R\right)&:=& \sup_{u\in W_r,\,\left\|
u\right\| =1 }\int_{\mathbb{R}^{N}\setminus B_{R}}K\left( \left| x\right|
\right) \left| u\right| ^{q}dx.  \label{S_i :=}
\end{eqnarray}
Clearly $\mathcal{S}_{0}\left( q,\cdot \right) $ is nondecreasing, $\mathcal{%
S}_{\infty }\left( q,\cdot \right) $ is nonincreasing and both of them can
be infinite at some $R$.

Our first result concerns the embedding properties of $W_{r}$ into the sum
space 
\[
L_{K}^{q_{1}}+L_{K}^{q_{2}}:=\left\{ u_{1}+u_{2}:u_{1}\in
L_{K}^{q_{1}}\left( \mathbb{R}^{N}\right) ,\,u_{2}\in L_{K}^{q_{2}}\left( \mathbb{R%
}^{N}\right) \right\} ,\quad 1<q_{i}<\infty . 
\]
We recall from \cite{BPR} that such a space can be characterized as the set
of measurable mappings $u:\mathbb{R}^{N}\rightarrow \mathbb{R}$ for which there
exists a measurable set $E\subseteq \mathbb{R}^{N}$ such that $u\in
L_{K}^{q_{1}}\left( E\right) \cap L_{K}^{q_{2}}\left( E^{c}\right) $. It is
a Banach space with respect to the norm 
\[
\left\| u\right\| _{L_{K}^{q_{1}}+L_{K}^{q_{2}}}:=\inf_{u_{1}+u_{2}=u}\max
\left\{ \left\| u_{1}\right\| _{L_{K}^{q_{1}}(\mathbb{R}^{N})},\left\|
u_{2}\right\| _{L_{K}^{q_{2}}(\mathbb{R}^{N})}\right\} 
\]
and the continuous embedding $L_{K}^{q}\hookrightarrow
L_{K}^{q_{1}}+L_{K}^{q_{2}}$ holds for all $q\in \left[ \min \left\{
q_{1},q_{2}\right\} ,\max \left\{ q_{1},q_{2}\right\} \right] $. The
assumptions of our result are quite general, sometimes also sharp (see claim
(iii)), but not so easy to check, so that the next results will be devoted
to provide more handy conditions ensuring such general assumptions.

\begin{thm}
\label{THM(cpt)} Let $1<p<N$, let $V$, $K$ be as in $\left( \mathbf{V}%
\right) $, $\left( \mathbf{K}\right) $ and let $q_{1},q_{2}>1$.

\begin{itemize}
\item[(i)]  If 
\begin{equation}
\mathcal{S}_{0}\left( q_{1},R_{1}\right) <\infty \quad \text{and}\quad 
\mathcal{S}_{\infty }\left( q_{2},R_{2}\right) <\infty \quad \text{for some }%
R_{1},R_{2}>0,  
\tag*{$\left( {\cal S}_{q_{1},q_{2}}^{\prime }\right) $}
\end{equation}
then $W_{r}$ is continuously embedded into $L_{K}^{q_{1}}(\mathbb{R}%
^{N})+L_{K}^{q_{2}}(\mathbb{R}^{N})$.

\item[(ii)]  If 
\begin{equation}
\lim_{R\rightarrow 0^{+}}\mathcal{S}_{0}\left( q_{1},R\right)
=\lim_{R\rightarrow +\infty }\mathcal{S}_{\infty }\left( q_{2},R\right) =0, 
\tag*{$\left( {\cal S}_{q_{1},q_{2}}^{\prime \prime }\right) $}
\end{equation}
then $W_{r}$ is compactly embedded into $L_{K}^{q_{1}}(\mathbb{R}%
^{N})+L_{K}^{q_{2}}(\mathbb{R}^{N})$.

\item[(iii)]  If $K\left( \left| \cdot \right| \right) \in L^{1}(B_{1})$ and 
$q_{1}\leq q_{2}$, then conditions $\left( \mathcal{S}_{q_{1},q_{2}}^{\prime
}\right) $ and $\left( \mathcal{S}_{q_{1},q_{2}}^{\prime \prime }\right) $
are also necessary to the above embeddings.
\end{itemize}
\end{thm}

Observe that, of course, $(\mathcal{S}_{q_{1},q_{2}}^{\prime \prime })$
implies $(\mathcal{S}_{q_{1},q_{2}}^{\prime })$. Moreover, these assumptions
can hold with $q_{1}=q_{2}=q$ and therefore Theorem \ref{THM(cpt)} also
concerns the embedding properties of $W_r$ into $L_{K}^{q}$, $1<q<\infty $.
\smallskip

We now look for explicit conditions on $V$ and $K$ implying $(\mathcal{S}%
_{q_{1},q_{2}}^{\prime \prime })$ for some $q_{1}$ and $q_{2}$. More
precisely, we will ensure $(\mathcal{S}_{q_{1},q_{2}}^{\prime \prime })$
through a more stringent condition involving the following functions of $R>0$
and $q>1$: 
\begin{eqnarray}
\mathcal{R}_{0}\left( q,R\right) := &&\sup_{u\in W_{r},\,h\in W,\,\left\|
u\right\| =\left\| h\right\| =1}\,\int_{B_{R}}K\left( \left| x\right|
\right) \left| u\right| ^{q-1}\left| h\right| dx,  \label{N_o} \\
\mathcal{R}_{\infty }\left( q,R\right) := &&\sup_{u\in W_{r},\,h\in
W,\,\left\| u\right\| =\left\| h\right\| =1}\,\int_{\mathbb{R}^{N}\setminus
B_{R}}K\left( \left| x\right| \right) \left| u\right| ^{q-1}\left| h\right|
dx.  \label{N_i}
\end{eqnarray}
Note that $\mathcal{R}_{0}\left( q,\cdot \right) $ is nondecreasing, $%
\mathcal{R}_{\infty }\left( q,\cdot \right) $ is nonincreasing and both can
be infinite at some $R$. Moreover, for every $\left( q,R\right) $ one has $%
\mathcal{S}_{0}\left( q,R\right) \leq \mathcal{R}_{0}\left( q,R\right) $ and 
$\mathcal{S}_{\infty }\left( q,R\right) \leq \mathcal{R}_{\infty }\left(
q,R\right) $, so that $(\mathcal{S}_{q_{1},q_{2}}^{\prime \prime })$ is a
consequence of the following, stronger condition: 
\begin{equation}
\lim_{R\rightarrow 0^{+}}\mathcal{R}_{0}\left( q_{1},R\right)
=\lim_{R\rightarrow +\infty }\mathcal{R}_{\infty }\left( q_{2},R\right) =0. 
\tag*{$\left( {\cal R}_{q_{1},q_{2}}^{\prime \prime }\right) $}
\end{equation}
In Theorems \ref{THM0} and \ref{THM3} we will find ranges of exponents $%
q_{1} $ such that $\lim_{R\rightarrow 0^{+}}\mathcal{R}_{0}\left(
q_{1},R\right) $ $=0$. In Theorems \ref{THM1} and \ref{THM2} we will do the
same for exponents $q_{2}$ such that $\lim_{R\rightarrow +\infty }\mathcal{R}%
_{\infty }\left( q_{2},R\right) =0$. Condition $(\mathcal{R}%
_{q_{1},q_{2}}^{\prime \prime })$ then follows by joining Theorem \ref{THM0}
or \ref{THM3} with Theorem \ref{THM1} or \ref{THM2}.\smallskip

For $\alpha \in \mathbb{R}$ and $\beta \in \left[ 0,1\right] $, define two
functions $\alpha ^{*}\left( \beta \right) $ and $q^{*}\left( \alpha ,\beta
\right) $ by setting 
\[
\alpha ^{*}\left( \beta \right) :=\max \left\{ p\beta -1-\frac{p-1}{p}N
,-\left( 1-\beta \right) N\right\} =\left\{ 
\begin{array}{ll}
p\beta -1-\frac{p-1}{p}N \quad \smallskip & \text{if }0\leq \beta \leq \frac{%
1}{p} \\ 
-\left( 1-\beta \right) N & \text{if }\frac{1}{p}\leq \beta \leq 1
\end{array}
\right. 
\]
and 
\[
q^{*}\left( \alpha ,\beta \right) :=p\frac{\alpha -p\beta +N}{N-p}. 
\]
Note that $\alpha ^{*}\left( \beta \right) \leq 0$ and $\alpha ^{*}\left(
\beta \right) =0$ if and only if $\beta =1$.

The first two Theorems \ref{THM0} and \ref{THM1} only rely on a power type
estimate of the relative growth of the potentials and do not require any
other separate assumption on $V$ and $K$ than $\left( \mathbf{V}\right) $
and $\left( \mathbf{K}\right) $, including the case $V=0$ (see Remark \ref
{RMK: suff12}.\ref{RMK: suff12-V^0}).

\begin{thm}
\label{THM0} Let $1<p<N$ and let $V$, $K$ be as in $\left( \mathbf{V}\right) 
$, $\left( \mathbf{K}\right) $. Assume that there exists $R_{1}>0$ such that 
$V\left( r\right) <+\infty $ almost everywhere in $(0,R_{1})$ and 
\begin{equation}
\esssup_{r\in \left( 0,R_{1}\right) }\frac{K\left( r\right) }{%
r^{\alpha _{0}}V\left( r\right) ^{\beta _{0}}}<+\infty \quad \text{for some }%
0\leq \beta _{0}\leq 1\text{~and }\alpha _{0}>\alpha ^{*}\left( \beta
_{0}\right) .  \label{esssup in 0}
\end{equation}
Then $\displaystyle \lim_{R\rightarrow 0^{+}}\mathcal{R}_{0}\left(
q_{1},R\right) =0$ for every $q_{1}\in \mathbb{R}$ such that 
\begin{equation}
\max \left\{ 1,p\beta _{0}\right\} <q_{1}<q^{*}\left( \alpha _{0},\beta
_{0}\right) .  \label{th1}
\end{equation}
\end{thm}

\begin{thm}
\label{THM1} Let $1<p<N$ and let $V$, $K$ be as in $\left( \mathbf{V}\right) 
$, $\left( \mathbf{K}\right) $. Assume that there exists $R_{2}>0$ such that 
$V\left( r\right) <+\infty $ for almost every $r>R_{2}$ and 
\begin{equation}
\esssup_{r>R_{2}}\frac{K\left( r\right) }{r^{\alpha _{\infty
}}V\left( r\right) ^{\beta _{\infty }}}<+\infty \quad \text{for some }0\leq
\beta _{\infty }\leq 1\text{~and }\alpha _{\infty }\in \mathbb{R}.
\label{esssup all'inf}
\end{equation}
Then $\displaystyle \lim_{R\rightarrow +\infty }\mathcal{R}_{\infty }\left(
q_{2},R\right) =0$ for every $q_{2}\in \mathbb{R}$ such that 
\begin{equation}
q_{2}>\max \left\{ 1,p\beta _{\infty },q^{*}\left( \alpha _{\infty },\beta
_{\infty }\right) \right\} .  \label{th2}
\end{equation}
\end{thm}

We observe explicitly that for every $\left( \alpha ,\beta \right) \in %
\mathbb{R}\times \left[ 0,1\right] $ one has 
\[
\max \left\{ 1,p\beta ,q^{*}\left( \alpha ,\beta \right) \right\} =\left\{ 
\begin{array}{ll}
q^{*}\left( \alpha ,\beta \right) \quad & \text{if }\alpha \geq \alpha
^{*}\left( \beta \right) \smallskip \\ 
\max \left\{ 1,p\beta \right\} & \text{if }\alpha \leq \alpha ^{*}\left(
\beta \right)
\end{array}
\right. . 
\]

\begin{rem}
\label{RMK: suff12}\quad

\begin{enumerate}
\item  \label{RMK: suff12-V^0}We mean $V\left( r\right) ^{0}=1$ for every $r$
(even if $V\left( r\right) =0$). In particular, if $V\left( r\right) =0$ for
almost every $r>R_{2}$, then Theorem \ref{THM1} can be applied with $\beta
_{\infty }=0$ and assumption (\ref{esssup all'inf}) means 
\[
\esssup_{r>R_{2}}\frac{K\left( r\right) }{r^{\alpha _{\infty }}}%
<+\infty \quad \text{for some }\alpha _{\infty }\in \mathbb{R}. 
\]
Similarly for Theorem \ref{THM0} and assumption (\ref{esssup in 0}), if $%
V\left( r\right) =0$ for almost every $r\in \left( 0,R_{1}\right) $.

\item  \label{RMK: suff12-no hp}The inequality $\max \left\{ 1,p\beta
_{0}\right\} <q^{*}\left( \alpha _{0},\beta _{0}\right) $ is equivalent to $%
\alpha _{0}>\alpha ^{*}\left( \beta _{0}\right) $. Then, in (\ref{th1}),
such inequality is automatically true and does not ask for further
conditions on $\alpha _{0}$ and $\beta _{0}$.

\item  \label{RMK: suff12-Vbdd}The assumptions of Theorems \ref{THM0} and 
\ref{THM1} may hold for different pairs $\left( \alpha _{0},\beta
_{0}\right) $, $\left( \alpha _{\infty },\beta _{\infty }\right) $. In this
case, of course, one chooses them in order to get the ranges for $%
q_{1},q_{2} $ as large as possible. For instance, if $V$ is essentially
bounded in a neighbourhood of 0 and condition (\ref{esssup in 0}) holds true
for a pair $\left( \alpha _{0},\beta _{0}\right) $, then (\ref{esssup in 0})
also holds for all pairs $\left( \alpha _{0}^{\prime },\beta _{0}^{\prime
}\right) $ such that $\alpha _{0}^{\prime }<\alpha _{0}$ and $\beta
_{0}^{\prime }<\beta _{0}$. Therefore, since $\max \left\{ 1,p\beta \right\} 
$ is nondecreasing in $\beta $ and $q^{*}\left( \alpha ,\beta \right) $ is
increasing in $\alpha $ and decreasing in $\beta $, it is convenient to
choose $\beta _{0}=0$ and the best interval where one can take $q_{1}$ is $%
1<q_{1}<q^{*}\left( \overline{\alpha },0\right) $ with $\overline{\alpha }%
:=\sup \{\alpha _{0}:\esssup_{r\in \left( 0,R_{1}\right) }K\left(
r\right) /r^{\alpha _{0}}<+\infty \}$ (we mean $q^{*}\left( +\infty
,0\right) =+\infty $).
\end{enumerate}
\end{rem}

For any $\alpha \in \mathbb{R}$, $\beta \leq 1$ and $\gamma \in \mathbb{R}$,
define 
\begin{equation}
q_{*}\left( \alpha ,\beta ,\gamma \right) :=p\frac{\alpha -\gamma \beta +N}{%
N-\gamma }\quad \text{and}\quad q_{**}\left( \alpha ,\beta ,\gamma \right)
:= p\frac{p\alpha +\left( 1-p\beta \right) \gamma +p\left( N-1\right) }{%
p\left( N-1\right) -\gamma (p-1)}.  \label{q** :=}
\end{equation}
Of course $q_{*}$ and $q_{**}$ are undefined if $\gamma =N$ and $\gamma = 
\frac{p}{p-1}\left( N-1\right) $, respectively.

The next Theorems \ref{THM2} and \ref{THM3} improve the results of Theorems 
\ref{THM0} and \ref{THM1} by exploiting further informations on the growth
of $V$ (see Remarks \ref{RMK: Hardy 1}.\ref{RMK: Hardy 1-improve} and \ref
{RMK: Hardy 2}.\ref{RMK: Hardy 2-improve}).

\begin{thm}
\label{THM2} Let $1<p<N$ and let $V$, $K$ be as in $\left( \mathbf{V}\right) 
$, $\left( \mathbf{K}\right) $. Assume that there exists $R_{2}>0$ such that 
$V\left( r\right) <+\infty $ for almost every $r>R_{2}$ and 
\begin{equation}
\esssup_{r>R_{2}}\frac{K\left( r\right) }{r^{\alpha _{\infty
}}V\left( r\right) ^{\beta _{\infty }}}<+\infty \quad \text{for some }0\leq
\beta _{\infty }\leq 1\text{~and }\alpha _{\infty }\in \mathbb{R}
\label{hp all'inf}
\end{equation}
and 
\begin{equation}
\essinf_{r>R_{2}}r^{\gamma _{\infty }}V\left( r\right) >0\quad 
\text{for some }\gamma _{\infty }\leq p.  \label{stima all'inf}
\end{equation}
Then $\displaystyle \lim_{R\rightarrow +\infty }\mathcal{R}_{\infty }\left(
q_{2},R\right) =0$ for every $q_{2}\in \mathbb{R}$ such that 
\begin{equation}
q_{2}>\max \left\{ 1,p\beta _{\infty },q_{*},q_{**}\right\} ,  \label{th3}
\end{equation}
where $q_{*}=q_{*}\left( \alpha _{\infty },\beta _{\infty },\gamma _{\infty
}\right) $ and $q_{**}=q_{**}\left( \alpha _{\infty },\beta _{\infty
},\gamma _{\infty }\right) .$
\end{thm}

For future convenience, we define three functions $\alpha _{1}:=\alpha
_{1}\left( \beta ,\gamma \right) $, $\alpha _{2}:=\alpha _{2}\left( \beta
\right) $ and $\alpha _{3}:=\alpha _{3}\left( \beta ,\gamma \right) $ by
setting 
\begin{equation}
\alpha _{1}:=-\left( 1-\beta \right) \gamma ,\quad \alpha _{2}:=-\left(
1-\beta \right) N,\quad \alpha _{3}:=-\frac{(p-1) N+\left( 1-p\beta \right)
\gamma }{p}.  \label{alpha_i :=}
\end{equation}
Then an explicit description of $\max \left\{ 1,p\beta ,q_{*},q_{**}\right\} 
$ is the following: for every $\left( \alpha ,\beta ,\gamma \right) \in %
\mathbb{R}\times \left( -\infty ,1\right] \times \left( -\infty ,N\right) $
we have 
\begin{equation}
\max \left\{ 1,p\beta ,q_{*},q_{**}\right\} =\left\{ 
\begin{array}{ll}
q_{**}\left( \alpha ,\beta ,\gamma \right) \quad & \text{if }\alpha \geq
\alpha _{1}\smallskip \\ 
q_{*}\left( \alpha ,\beta ,\gamma \right) & \text{if }\max \left\{ \alpha
_{2},\alpha _{3}\right\} \leq \alpha \leq \alpha _{1}\smallskip \\ 
\max \left\{ 1,p\beta \right\} & \text{if }\alpha \leq \max \left\{ \alpha
_{2},\alpha _{3}\right\}
\end{array}
\right. ,  \label{descrizioneThm2}
\end{equation}
where $\max \left\{ \alpha _{2},\alpha _{3}\right\} <\alpha _{1}$ for every $%
\beta <1$ and $\max \left\{ \alpha _{2},\alpha _{3}\right\} =\alpha _{1}=0$
if $\beta =1$.

\begin{rem}
\label{RMK: Hardy 1}\quad

\begin{enumerate}
\item  \label{RMK: Hardy 1-B<0}The proof of Theorem \ref{THM2} does not
require $\beta _{\infty }\geq 0$, but this condition is not a restriction of
generality in stating the theorem. Indeed, under assumption (\ref{stima
all'inf}), if (\ref{hp all'inf}) holds with $\beta _{\infty }<0$, then it
also holds with $\alpha _{\infty }$ and $\beta _{\infty }$ replaced by $%
\alpha _{\infty }-\beta _{\infty }\gamma _{\infty }$ and $0$ respectively,
and this does not change the thesis (\ref{th3}), because $q_{*}\left( \alpha
_{\infty }-\beta _{\infty }\gamma _{\infty },0,\gamma _{\infty }\right)
=q_{*}\left( \alpha _{\infty },\beta _{\infty },\gamma _{\infty }\right) $
and $q_{**}\left( \alpha _{\infty }-\beta _{\infty }\gamma _{\infty
},0,\gamma _{\infty }\right) =q_{**}\left( \alpha _{\infty },\beta _{\infty
},\gamma _{\infty }\right) $.

\item  \label{RMK: Hardy 1-improve}

Denote $q^{*}=q^{*}\left( \alpha _{\infty },\beta _{\infty }\right) $ for
brevity. If $\gamma _{\infty }<p$, then one has 
\[
\max \left\{ 1,p\beta _{\infty },q^{*}\right\} =\left\{ 
\begin{array}{ll}
\max \left\{ 1,p\beta _{\infty }\right\} =\max \left\{ 1,p\beta _{\infty
},q_{*},q_{**}\right\} \quad \smallskip & \text{if }\alpha _{\infty }\leq
\alpha ^{*}\left( \beta _{\infty }\right) \\ 
q^{*}>\max \left\{ 1,p\beta _{\infty },q_{*},q_{**}\right\} & \text{if }%
\alpha _{\infty }>\alpha ^{*}\left( \beta _{\infty }\right)
\end{array}
\right. , 
\]
so that, under assumption (\ref{stima all'inf}), Theorem \ref{THM2} improves
Theorem \ref{THM1}. Otherwise, if $\gamma _{\infty }=p$, we have $%
q_{*}=q_{**}=q^{*}$ and Theorems \ref{THM2} and \ref{THM1} give the same
result. This is not surprising, since, by Hardy inequality, the space $W$
coincides with $D^{1,p}(\mathbb{R}^{N})$ if $V\left( r\right) =r^{-p}$ and
thus, for $\gamma _{\infty }=p$, we cannot expect a better result than the
one of Theorem \ref{THM1}, which covers the case of $V=0$, i.e., of $D^{1,p}(%
\mathbb{R}^{N})$.

\item  \label{RMK: Hardy 1-best gamma}Description (\ref{descrizioneThm2})
shows that $q_{*}$ and $q_{**}$ are not relevant in inequality (\ref{th3})
if $\alpha _{\infty }\leq \alpha _{2}\left( \beta _{\infty }\right) $. On
the other hand, if $\alpha _{\infty }>\alpha _{2}\left( \beta _{\infty
}\right) $, both $q_{*}$ and $q_{**}$ turn out to be increasing in $\gamma $
and hence it is convenient to apply Theorem \ref{THM2} with the smallest $%
\gamma _{\infty }$ for which (\ref{stima all'inf}) holds. This is consistent
with the fact that, if (\ref{stima all'inf}) holds with $\gamma _{\infty }$,
then it also holds with every $\gamma _{\infty }^{\prime }$ such that $%
\gamma _{\infty }\leq \gamma _{\infty }^{\prime }\leq p$.
\end{enumerate}
\end{rem}

In order to state our last result, we introduce, by the following
definitions, an open region $\mathcal{A}_{\beta ,\gamma }$ of the $\alpha q$%
-plane, depending on $\beta\in[0,1]$ and $\gamma \geq p$. Recall the
definitions (\ref{q** :=}) of the functions $q_{*}=q_{*}\left( \alpha ,\beta
,\gamma \right) $ and $q_{**}=q_{**}\left( \alpha ,\beta ,\gamma \right) $.
We set 
\begin{equation}
\begin{array}{ll}
\mathcal{A}_{\beta ,\gamma }:=\left\{ \left( \alpha ,q\right) :\max \left\{
1,p\beta \right\} <q<\min \left\{ q_{*},q_{**}\right\} \right\} \quad
\smallskip & \text{if }p\leq\gamma <N, \\ 
\mathcal{A}_{\beta ,\gamma }:=\left\{ \left( \alpha ,q\right) :\max \left\{
1,p\beta \right\} <q<q_{**},\,\alpha >-\left( 1-\beta \right) N\right\}
\quad \smallskip & \text{if }\gamma =N, \\ 
\mathcal{A}_{\beta ,\gamma }:=\left\{ \left( \alpha ,q\right) :\max \left\{
1,p\beta ,q_{*}\right\} <q<q_{**}\right\} \smallskip & \text{if }N<\gamma <%
\frac{p}{p-1}(N-1), \\ 
\mathcal{A}_{\beta ,\gamma }:=\left\{ \left( \alpha ,q\right) :\max \left\{
1,p\beta ,q_{*}\right\} <q,\,\alpha >-\left( 1-\beta \right) \gamma \right\}
\smallskip & \text{if }\gamma =\frac{p}{p-1}(N-1), \\ 
\mathcal{A}_{\beta ,\gamma }:=\left\{ \left( \alpha ,q\right) :\max \left\{
1,p\beta ,q_{*},q_{**}\right\} <q\right\} & \text{if }\gamma >\frac{p}{p-1}%
(N-1).
\end{array}
\label{A:=}
\end{equation}
Notice that $\frac{p}{p-1}(N-1)>N$ because $p<N$. For more clarity, $%
\mathcal{A}_{\beta ,\gamma }$ is sketched in the following five pictures,
according to the five cases above. Recall the definitions (\ref{alpha_i :=})
of the functions $\alpha _{1}=\alpha _{1}\left( \beta ,\gamma \right) $, $%
\alpha _{2}=\alpha _{2}\left( \beta \right) $ and $\alpha _{3}=\alpha
_{3}\left( \beta ,\gamma \right) $.\bigskip


\noindent
\begin{tabular}[t]{l}
\begin{tabular}{l}
\textbf{Fig.1}:\ \ $\mathcal{A}_{\beta ,\gamma }$ for 
$p\leq \gamma <N$.\smallskip\\
$\bullet $ If $\gamma =p$, the two straight \\
lines above are the same.\smallskip\\
$\bullet $ If $\beta <1$ we have\\
$\max \left\{ \alpha _{2},\alpha _{3}\right\} <\alpha _{1}<0.$\\
If $\beta =1$ we have\\
$\alpha _{3}<\alpha _{2}=\alpha _{1}=0$\\
and $\mathcal{A}_{1,\gamma }$ reduces to the angle\\
$p<q<q_{**}$.
\end{tabular}
\begin{tabular}{l}
\includegraphics[width=3.7in]{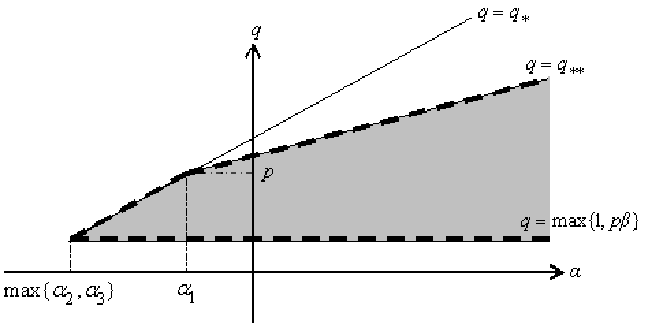}
\end{tabular}
\end{tabular}
\bigskip

\noindent
\begin{tabular}[t]{l}
\begin{tabular}{l}
\textbf{Fig.2}:\ \ $\mathcal{A}_{\beta ,\gamma }$ 
$\gamma =N$.\smallskip\\
$\bullet $ If $\beta <1$ we have\\
$\alpha _{1}=\alpha _{2}=\alpha _{3}<0.$\\
If $\beta =1$ we have\\
$\alpha _{1}=\alpha _{2}=\alpha _{3}=0$\\
and $\mathcal{A}_{1,\gamma }$ reduces to the angle\\
$p<q<q_{**}$.
\end{tabular}
\begin{tabular}{l}
\includegraphics[width=3.7in]{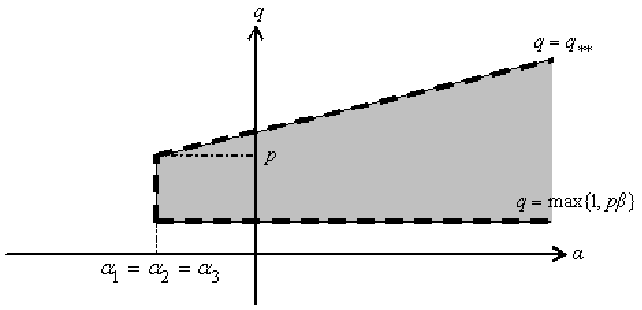}
\end{tabular}
\end{tabular}
\bigskip

\noindent
\begin{tabular}[t]{l}
\begin{tabular}{l}
\textbf{Fig.3}:\ \ $\mathcal{A}_{\beta ,\gamma }$ for \smallskip\\
$N<\gamma<\frac{p}{p-1}(N-1)$.\smallskip\\
$\bullet $ If $\beta <1$ we have\\
$\alpha _{1}<\min \left\{ \alpha _{2},\alpha _{3}\right\} <0.$\\
If $\beta =1$ we have\\
$0=\alpha _{1}=\alpha_{2}<\alpha_{3}$\\
and $\mathcal{A}_{1,\gamma }$ reduces to the angle\\
$p<q<q_{**}$.
\end{tabular}
\begin{tabular}{l}
\includegraphics[width=3.7in]{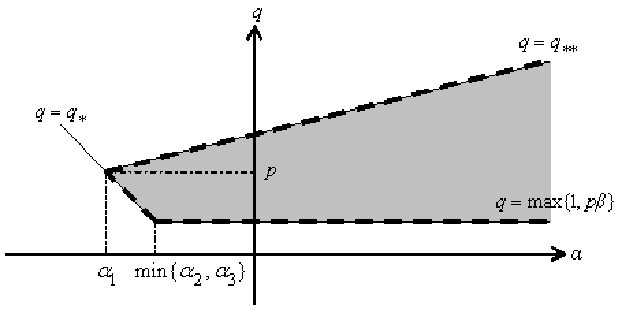}
\end{tabular}
\end{tabular}
\bigskip

\noindent
\begin{tabular}[t]{l}
\begin{tabular}{l}
\textbf{Fig.4}:\ \ $\mathcal{A}_{\beta ,\gamma }$ for 
$\gamma=\frac{p}{p-1}(N-1)$.\smallskip\\ 
$\bullet $ If $\beta <1$ we have \\ 
$\alpha _{1}<\min \left\{ \alpha _{2},\alpha _{3}\right\} <0.$ \\ 
If $\beta =1$ we have \\ 
$0=\alpha _{1}=\alpha_{2}<\alpha_{3}$ \\ 
and $\mathcal{A}_{1,\gamma }$ reduces to the angle\\ 
$\alpha >0,\,q>p$.
\end{tabular}
\begin{tabular}{l}
\includegraphics[width=3.4in]{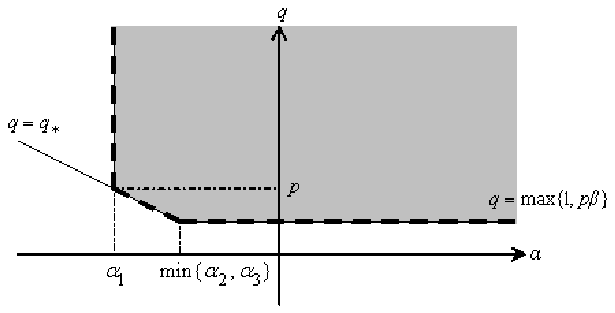}
\end{tabular}
\end{tabular}
\bigskip

\noindent
\begin{tabular}[t]{l}
\begin{tabular}{l}
\textbf{Fig.5}:\ \ $\mathcal{A}_{\beta ,\gamma }$ for 
$\gamma>\frac{p}{p-1}(N-1)$.\smallskip\\
$\bullet $ If $\beta <1$ we have\\
$\alpha _{1}<\min \left\{ \alpha _{2},\alpha _{3}\right\} <0.$\\
If $\beta =1$ we have\\
$0=\alpha _{1}=\alpha_{2}<\alpha_{3}$\\
and $\mathcal{A}_{1,\gamma }$ reduces to the angle\\
$q>\max \left\{p,q_{**}\right\} $.
\end{tabular}
\begin{tabular}{l}
\includegraphics[width=3.4in]{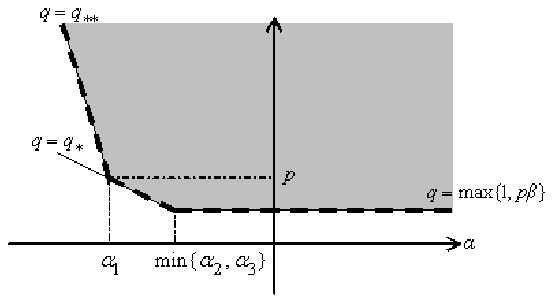}
\end{tabular}
\end{tabular}
\bigskip

\begin{thm}
\label{THM3}Let $N\geq 3$ and let $V$, $K$ be as in $\left( \mathbf{V}%
\right) $, $\left( \mathbf{K}\right) $. 
Assume that there exists $R_{1}>0$ such that $V\left( r\right) <+\infty $
almost everywhere in $(0,R_{1})$ and 
\begin{equation}
\esssup_{r\in \left( 0,R_{1}\right) }\frac{K\left( r\right) }{%
r^{\alpha _{0}}V\left( r\right) ^{\beta _{0}}}<+\infty \quad \text{for some }%
0\leq \beta _{0}\leq 1\text{~and }\alpha _{0}\in \mathbb{R}  \label{hp in 0}
\end{equation}
and 
\begin{equation}
\essinf_{r\in \left( 0,R_{1}\right) }r^{\gamma _{0}}V\left(
r\right) >0\quad \text{for some }\gamma _{0}\geq p.  \label{stima in 0}
\end{equation}
Then $\displaystyle \lim_{R\rightarrow 0^{+}}\mathcal{R}_{0}\left(
q_{1},R\right) =0$ for every $q_{1}\in \mathbb{R}$ such that 
\begin{equation}
\left( \alpha _{0},q_{1}\right) \in \mathcal{A}_{\beta _{0},\gamma _{0}}.
\label{th4}
\end{equation}
\end{thm}

\begin{rem}
\label{RMK: Hardy 2}\quad

\begin{enumerate}
\item  Condition (\ref{th4}) also asks for a lower bound on $\alpha _{0}$,
except for the case $\gamma _{0}>\frac{p}{p-1}(N-1)$, as it is clear from
Figures 1-5.

\item  The proof of Theorem \ref{THM3} does not require $\beta _{0}\geq 0$,
but this is not a restriction of generality in stating the theorem (cf.
Remark \ref{RMK: Hardy 1}.\ref{RMK: Hardy 1-B<0}). Indeed, under assumption (%
\ref{stima in 0}), if (\ref{hp in 0}) holds with $\beta _{0}<0$, then it
also holds with $\alpha _{0}$ and $\beta _{0}$ replaced by $\alpha
_{0}-\beta _{0}\gamma _{0}$ and $0$ respectively, and one has that $\left(
\alpha _{0},q_{1}\right) \in \mathcal{A}_{\beta _{0},\gamma _{0}}$ if and
only if $\left( \alpha _{0}-\beta _{0}\gamma _{0},q_{1}\right) \in \mathcal{A%
}_{0,\gamma _{0}}$.

\item  \label{RMK: Hardy 2-improve}If (\ref{stima in 0}) holds with $\gamma
_{0}>p$, then Theorem \ref{THM3} improves Theorem \ref{THM0}. Otherwise, if $%
\gamma _{0}=p$, then one has $\max \left\{ \alpha _{2},\alpha _{3}\right\}
=\alpha ^{*}\left( \beta _{0}\right) $ and $\left( \alpha _{0},q_{1}\right)
\in \mathcal{A}_{\beta _{0},\gamma _{0}}$ is equivalent to $\max \left\{
1,p\beta _{0}\right\} <q_{1}<q^{*}\left( \alpha _{0},\beta _{0}\right) $,
i.e., Theorems \ref{THM3} and \ref{THM0} give the same result, which is
consistent with Hardy inequality (cf. Remark \ref{RMK: Hardy 1}.\ref{RMK:
Hardy 1-improve}).

\item  \label{RMK: Hardy 2-best gamma}Given $\beta \leq 1$, one can check
that $\mathcal{A}_{\beta ,\gamma _{1}}\subseteq \mathcal{A}_{\beta ,\gamma
_{2}}$ for every $p\leq \gamma _{1}<\gamma _{2}$, so that, in applying
Theorem \ref{THM3}, it is convenient to choose the largest $\gamma _{0}$ for
which (\ref{stima in 0}) holds. This is consistent with the fact that, if (%
\ref{stima in 0}) holds with $\gamma _{0}$, then it also holds with every $%
\gamma _{0}^{\prime }$ such that $p\leq \gamma _{0}^{\prime }\leq \gamma
_{0} $.
\end{enumerate}
\end{rem}

\begin{rem}
If $p=2$, the above compactness theorems exactly reduces to the ones of \cite
{BGR_I}, except for the fact that there we required assumption $\left( 
\mathbf{K}\right) $ with $s>2N/(N+2)$ instead of $s>1$. In this respect, the
result we present here are improvements of the ones of \cite{BGR_I} also for 
$p=2$.
\end{rem}


\section{Examples \label{SEC: examples}}

In this section we give some examples of application of our compactness
results, which might clarify how to use them in concrete cases. We also
compare them with the most recent and general related results \cite
{Su-Wang-Will-p,SuTian12}, which unify and extend the previous literature.
Essentially, the spirit of the results of \cite{Su-Wang-Will-p,SuTian12} is
the following: assuming that $V,K$ are continuous and satisfy power type
estimates of the form: 
\begin{equation}
\liminf_{r\rightarrow 0^{+}}\frac{V\left( r\right) }{r^{a_{0}}}>0,\quad
\liminf_{r\rightarrow +\infty }\frac{V\left( r\right) }{r^{a}}>0,\quad
\limsup_{r\rightarrow 0^{+}}\frac{K\left( r\right) }{r^{b_{0}}}<\infty
,\quad \limsup_{r\rightarrow +\infty }\frac{K\left( r\right) }{r^{b}}<\infty
,  \label{pow-estim}
\end{equation}
the authors find two limit exponents $\underline{q}=\underline{q}\left(
a,b\right) $ and $\overline{q}=\overline{q}\left( a_{0},b_{0}\right) $ such
that the embedding $W_{r}\hookrightarrow L_{K}^{q}$ is compact if $%
\underline{q}<q<\overline{q}$. The case with $q>p$ is studied in \cite
{Su-Wang-Will-p}, the one with $q<p$ in \cite{SuTian12}. The exponent $%
\underline{q}$ is always defined, while $\overline{q}$ exists provided that
suitable compatibility conditions between $a_{0}$ and $b_{0}$ occur.
Moreover, the condition $\underline{q}<q<\overline{q}$ also asks for $%
\underline{q}<\overline{q}$, which is a further assumption of compatibility
between the behaviours of the potentials at zero and at infinity.

In the following it will be always understood that $1<p<N$.

\begin{example}
\label{EX(SWW)}Consider the potentials 
\[
V\left( r\right) =\frac{1}{r^{a}},\quad K\left( r\right) =\frac{1}{r^{a-1}}%
,\quad a\leq p. 
\]
Since $V$ satisfies (\ref{stima all'inf}) with $\gamma _{\infty }=a$ (cf.
Remark \ref{RMK: Hardy 1}.\ref{RMK: Hardy 1-best gamma} for the best choice
of $\gamma _{\infty }$), we apply Theorems \ref{THM0} and \ref{THM2}, where
we choose $\beta _{0}=\beta _{\infty }=0$ and $\alpha _{0}=\alpha _{\infty
}=1-a$. Note that $a\leq p$ implies $\alpha _{0}=1-a>\alpha ^{*}\left(
0\right) $ and $q_{*}\left( \alpha _{\infty },0,a\right) \leq q_{**}\left(
\alpha _{\infty },0,a\right) $. Hence we get that $\left( \mathcal{R}%
_{q_{1},q_{2}}^{\prime \prime }\right) $ holds for every exponents $%
q_{1},q_{2}$ such that 
\begin{equation}
1<q_{1}<q^{*}=p\frac{N-a+1}{N-p},\quad q_{2}>q_{**}=p\frac{pN-a\left(
p-1\right) }{p\left( N-1\right) -a(p-1)}.  \label{q1q2}
\end{equation}
If $a<p$, then one has $q_{**}<q^{*}$ and therefore Theorem \ref{THM(cpt)}
gives the compact embedding 
\begin{equation}
W_{r}\hookrightarrow L_{K}^{q}\qquad \text{for}\quad q_{**}<q<q^{*}.
\label{ES(sww): p}
\end{equation}
If $a=p$, then $q_{**}=q^{*}$ and we get the compact embedding 
\[
W_{r}\hookrightarrow L_{K}^{q_{1}}+L_{K}^{q_{2}}\qquad \text{for}\quad
1<q_{1}<p+\frac{p}{N-p}<q_{2}. 
\]
Since $V$ and $K$ are power potentials, one can also apply the results of 
\cite{Su-Wang-Will-p}, which give two suitable limit exponents $\underline{q}
$ and $\overline{q}$ such that the embedding $W_{r}\hookrightarrow L_{K}^{q}$
is compact if $\underline{q}<q<\overline{q}$. These exponents $\underline{q}$
and $\overline{q}$ are exactly exponents $q_{**}$ and $q^{*}$ of (\ref{q1q2}%
) respectively, so that one obtains (\ref{ES(sww): p}) again provided that $%
a<p$ (which implies $\underline{q}<\overline{q}$). If $a=p$, instead, one
gets $\underline{q}=\overline{q}$ and no result is avaliable in \cite
{Su-Wang-Will-p}. The results of \cite{SuTian12} do not apply to $V$ and $K$%
, since the top and bottom exponents of \cite{SuTian12} turn out to be equal
to one another for every $a\leq p$.
\end{example}

The next Examples \ref{EX: BPR}, \ref{EX(nnP1)} and \ref{EX(nnP2)} concern
potentials for which no result is available in \cite{Su-Wang-Will-p,SuTian12}%
, since they do not satisfy (\ref{pow-estim}).

\begin{example}
\label{EX: BPR}Taking $V=0$, $K$ as in $\left( \mathbf{K}\right) $ and $%
\beta _{0}=\beta _{\infty }=0$, from Theorems \ref{THM0} and \ref{THM1} (see
also Remark \ref{RMK: suff12}.\ref{RMK: suff12-V^0}) we get that $\left( 
\mathcal{R}_{q_{1},q_{2}}^{\prime \prime }\right) $ holds for 
\begin{equation}
1<q_{1}<p\frac{\alpha _{0}+N}{N-p}\quad \text{and}\quad q_{2}>\max \left\{
1,p\frac{\alpha _{\infty }+N}{N-p}\right\} ,  \label{EX(BPR): th}
\end{equation}
provided that $\exists R_{1},R_{2}>0$ such that 
\begin{equation}
\esssup_{r>R_{2}}\frac{K\left( r\right) }{r^{\alpha _{\infty }}}%
<+\infty \quad \text{and}\quad \esssup_{r\in \left( 0,R_{1}\right) }%
\frac{K\left( r\right) }{r^{\alpha _{0}}}<+\infty \quad \text{with }\alpha
_{0}>-1-\frac{p-1}{p}N.  \label{EX(BPR): hp}
\end{equation}
Correspondingly, Theorem \ref{THM(cpt)} gives a compact embedding of $D_{%
\mathrm{rad}}^{1,p}(\mathbb{R}^{N})$ into $L_{K}^{q_{1}}+L_{K}^{q_{2}}$, which
was already proved in \cite[Theorem 4.1]{BPR} assuming $K\in L_{\mathrm{loc}%
}^{\infty }(\mathbb{R}_{+})$. Of course, according to (\ref{EX(BPR): th}), in (%
\ref{EX(BPR): hp}) it is convenient to choose $\alpha _{0}$ as large as
possible and $\alpha _{\infty }$ as small as possible. For instance, if $%
K\left( r\right) =r^{d}$ with $d>-1-N(p-1)/p$, we choose $\alpha _{0}=\alpha
_{\infty }=d$ and obtain the compact embedding 
\[
D_{\mathrm{rad}}^{1,2}(\mathbb{R}^{N})\hookrightarrow
L_{K}^{q_{1}}+L_{K}^{q_{2}}\qquad \text{for}\quad 1<q_{1}<p\frac{d+N}{N-p}%
<q_{2}. 
\]
Observe that, if (\ref{EX(BPR): hp}) holds for some $\alpha _{0}>\alpha
_{\infty }$, then we can take $q_{1}=q_{2}$ in (\ref{EX(BPR): th}) and get
the compact embedding 
\[
D_{\mathrm{rad}}^{1,2}(\mathbb{R}^{N})\hookrightarrow L_{K}^{q}\qquad \text{for}%
\quad \max \left\{ 1,p\frac{\alpha _{\infty }+N}{N-p}\right\} <q<p\frac{%
\alpha _{0}+N}{N-p}. 
\]
\end{example}

\begin{example}
\label{EX(nnP1)}Essentially the same result of Example \ref{EX: BPR} holds
if $V$ is not singular at the origin and, roughly speaking, decays at
infinity much faster than $K$ (or is compactly supported). The result
becomes different (and better) if $K$ decays at infinity similarly to $V$,
or much faster. For example, consider the potentials 
\[
V\left( r\right) =e^{-ar},\quad K_{1}\left( r\right) =r^{d},\quad
K_{2}\left( r\right) =r^{d}e^{-br},\quad a,b>0,\quad d>-1-\frac{p-1}{p}N. 
\]
Since $V$ does not satisfy (\ref{stima in 0}) or (\ref{stima all'inf}), we
use Theorems \ref{THM0} and \ref{THM1}. According to Remark \ref{RMK: suff12}%
.\ref{RMK: suff12-Vbdd}, both for $K=K_{1}$ and $K=K_{2}$, Theorem \ref{THM0}
leads to take $1<q_{1}<p(d+N)/(N-p)$. If $K=K_{1}$, the ratio in (\ref
{esssup all'inf}) is bounded only if $\beta _{\infty }=0$ and the best $%
\alpha _{\infty }$ we can take is $\alpha _{\infty }=d$, which yields $%
q_{2}>p(d+N)/(N-p)$. Then, via condition $\left( \mathcal{R}%
_{q_{1},q_{2}}^{\prime \prime }\right) $, Theorem \ref{THM(cpt)} gives the
compact embedding 
\[
W_{r}\hookrightarrow L_{K_{1}}^{q_{1}}+L_{K_{1}}^{q_{2}}\qquad \text{for}%
\quad 1<q_{1}<p\frac{d+N}{N-p}<q_{2}. 
\]
If $K=K_{2}$, instead, assumption (\ref{esssup all'inf}) holds with $\beta
_{\infty }=0$ and $\alpha _{\infty }\in \mathbb{R}$ arbitrary, so that we can
take $q_{2}>1$ arbitrary. Then, via condition $\left( \mathcal{R}%
_{q_{1},q_{2}}^{\prime \prime }\right) $ with $q_{2}=q_{1}$, Theorem \ref
{THM(cpt)} gives the compact embedding 
\[
W_{r}\hookrightarrow L_{K_{2}}^{q}\qquad \text{for}\quad 1<q<p\frac{d+N}{N-p}%
. 
\]
\end{example}

\begin{example}
\label{EX(nnP2)}Consider the potentials 
\[
V\left( r\right) =e^{\frac{1}{r}},\quad K\left( r\right) =e^{\frac{b}{r}%
},\quad 0<b\leq 1. 
\]
Since $V$ satisfies (\ref{stima in 0}), we apply Theorem \ref{THM(cpt)}
together with Theorems \ref{THM1} and \ref{THM3}. Assumption (\ref{esssup
all'inf}) holds for $\alpha _{\infty }\geq 0$ and $0\leq \beta _{\infty
}\leq 1$, so that the best choice for $\alpha _{\infty }$, which is $\alpha
_{\infty }=0$, gives 
\[
\max \left\{ 1,p\beta _{\infty },q^{*}\left( 0,\beta _{\infty }\right)
\right\} =p\frac{N-p\beta _{\infty }}{N-p}. 
\]
Then we take $\beta _{\infty }=1$, so that Theorem \ref{THM1} gives $q_{2}>p$%
. As to Theorem \ref{THM3}, hypothesis (\ref{stima in 0}) holds with $\gamma
_{0}\geq p$ arbitrary and therefore the most convenient choice is to assume $%
\gamma _{0}>(N-1)p/(p-1)$ (see Remark \ref{RMK: Hardy 2}.\ref{RMK: Hardy
2-best gamma}). On the other hand, we have 
\[
\frac{K\left( r\right) }{r^{\alpha _{0}}V\left( r\right) ^{\beta _{0}}}=%
\frac{e^{\frac{b-\beta _{0}}{r}}}{r^{\alpha _{0}}} 
\]
and thus hypothesis (\ref{hp in 0}) holds for some $\alpha _{0}\in \mathbb{R}$
if and only if $b\leq \beta _{0}\leq 1$. We now distinguish two cases. If $%
b<1$, we can take $\beta _{0}>b$ and thus (\ref{hp in 0}) holds for every $%
\alpha _{0}\in \mathbb{R}$, so that Theorem \ref{THM3} gives $q_{1}>\max
\left\{ 1,p\beta _{0}\right\} $ (see Fig.5), i.e., $q_{1}>\max \left\{
1,2b\right\} $. If $b=1$, then we need to take $\beta _{0}=1$ and thus (\ref
{hp in 0}) holds for $\alpha _{0}\leq 0$. Since $\gamma _{0}>(N-1)p/(p-1)$
implies 
\[
\mathcal{A}_{1,\gamma _{0}}=\left\{ \left( \alpha ,q\right) \in \mathbb{R}%
^{2}:q>\max \left\{ p,p-\frac{\alpha p^{2}}{\gamma _{0}\left( p-1\right)
-p\left( N-1\right) }\right\} \right\} , 
\]
the best choice for $\alpha _{0}\leq 0$ is $\alpha _{0}=0$ and we get that $%
\left( 0,q_{1}\right) \in \mathcal{A}_{1,\gamma _{0}}$ if and only if $%
q_{1}>p$. Hence Theorem \ref{THM3} gives $q_{1}>\max \left\{ 1,2b\right\} $
again. As a conclusion, observing that $0<b\leq 1$ implies $\max \left\{
1,pb\right\} \leq p$, we obtain condition $\left( \mathcal{R}_{q,q}^{\prime
\prime }\right) $ and the compact embedding $W_{r}\hookrightarrow L_{K}^{q}$
for $q>p$.

If we now modify $V$ by taking a compactly supported potential $V_{1}$ such
that $V_{1}\left( r\right) \sim V\left( r\right) $ as $r\rightarrow 0^{+}$,
everything works as above in applying Theorem \ref{THM3}, but now we need to
take $\beta _{\infty }=0$ and $\alpha _{\infty }\geq 0$ in Theorem \ref{THM1}%
. This gives 
\[
\max \left\{ 1,p\beta _{\infty },q^{*}\left( \alpha _{\infty },\beta
_{\infty }\right) \right\} =p\frac{\alpha _{\infty }+N}{N-p} 
\]
and thus, choosing $\alpha _{\infty }=0$, we get $\left( \mathcal{R}%
_{q,q}^{\prime \prime }\right) $ and the compact embedding $%
W_{r}\hookrightarrow L_{K}^{q}$ for $q>p^{*}$.

Similarly, if we modify $V$ by taking a potential $V_{2}$ such that $%
V_{2}\left( r\right) \sim V\left( r\right) $ as $r\rightarrow 0^{+}$ and $%
V_{2}\left( r\right) \sim r^{N}$ as $r\rightarrow +\infty $, Theorem \ref
{THM3} yields $q_{1}>\max \left\{ 1,pb\right\} $ as above and Theorem \ref
{THM1} gives $q_{2}>1$ (apply it for instance with $\alpha _{\infty }=-N/2$
and $\beta _{\infty }=1/2$), so that we get $\left( \mathcal{R}%
_{q,q}^{\prime \prime }\right) $ and the compact embedding $%
W_{r}\hookrightarrow L_{K}^{q}$ for $q>\max \left\{ 1,pb\right\} $.
\end{example}

The last example shows that our results also extend the ones of \cite
{Su-Wang-Will-p,SuTian12} for power-type potentials.

\begin{example}
\label{EX(ST)}Consider the potential 
\begin{equation}
V\left( r\right) =r^{a},\quad -\frac{p}{p-1}\left( N-1\right) <a<-N,
\label{EX1: V}
\end{equation}
and let $K$ be as in $\left( \mathbf{K}\right) $ and such that 
\begin{equation}
K\left( r\right) =O\left( r^{b_{0}}\right) _{r\rightarrow 0^{+}},\quad
K\left( r\right) =O\left( r^{b}\right) _{r\rightarrow +\infty },\quad
b_{0}>a,\quad b\in \mathbb{R}.  \label{EX1: K}
\end{equation}
Since $V$ satisfies (\ref{stima in 0}) with $\gamma _{0}=-a$ (cf. Remark \ref
{RMK: Hardy 2}.\ref{RMK: Hardy 2-best gamma} for the best choice of $\gamma
_{0}$), we apply Theorem \ref{THM(cpt)} together with Theorems \ref{THM1}
and \ref{THM3}, where assumptions (\ref{esssup all'inf}) and (\ref{hp in 0})
hold for $\alpha _{\infty }\geq b-a\beta _{\infty }$ and $\alpha _{0}\leq
b_{0}-a\beta _{0}$ with $0\leq \beta _{\infty }\leq 1$ and $\beta _{0}\leq 1$
arbitrary. Note that $N<\gamma _{0}<(N-1)p/(p-1)$. According to (\ref{th2})
and (\ref{th4}) (see in particular Fig.3), it is convenient to choose $%
\alpha _{\infty }$ as small as possible and $\alpha _{0}$ as large as
possible, so we take 
\begin{equation}
\alpha _{\infty }=b-a\beta _{\infty },\quad \alpha _{0}=b_{0}-a\beta _{0}.
\label{EX(SWW2):alpha}
\end{equation}
Then $q^{*}=q^{*}\left( \alpha _{\infty },\beta _{\infty }\right) $, $%
q_{*}=q_{*}\left( \alpha _{0},\beta _{0},-a\right) $ and $%
q_{**}=q_{**}\left( \alpha _{0},\beta _{0},-a\right) $ are given by 
\[
q^{*}=p\frac{N+b-\left( a+p\right) \beta _{\infty }}{N-p},\quad q_{*}=p\frac{%
N+b_{0}}{N+a}\quad \text{and}\quad q_{**}=p\frac{p(N-1)+pb_{0}-a}{p\left(
N-1\right) +a(p-1)}. 
\]
Since $a+p<p-N<0$, the exponent $q^{*}$ is increasing in $\beta _{\infty }$
and thus, according to (\ref{th2}) again, the best choice for $\beta
_{\infty }$ is $\beta _{\infty }=0$. This yields 
\begin{equation}
q_{2}>\max \left\{ 1,p\frac{N+b}{N-p}\right\} .  \label{EX(SWW2):q1}
\end{equation}
As to Theorem \ref{THM3}, we observe that, thanks to the choice of $\alpha
_{0}$, the exponents $q_{*}$ and $q_{**}$ are independent of $\beta _{0}$,
so that we can choose $\beta _{0}=0$ in order to get the region $\mathcal{A}%
_{\beta _{0},-a}$ as large as possible (cf. Fig.3 or the third definition in
(\ref{A:=})). Then we get $\alpha _{0}=b_{0}>a=\alpha _{1}$ (recall (\ref
{EX(SWW2):alpha}) and the definition (\ref{alpha_i :=}) of $\alpha _{1}$),
so that $\left( \alpha _{0},q_{1}\right) \in \mathcal{A}_{0,-a}$ if and only
if 
\begin{equation}
\max \left\{ 1,p\frac{N+b_{0}}{N+a}\right\} <q_{1}<p\frac{p(N-1)+pb_{0}-a}{%
p\left( N-1\right) +a(p-1)}.  \label{EX(SWW2):q2}
\end{equation}
As a conclusion, via condition $\left( \mathcal{R}_{q_{1},q_{2}}^{\prime
\prime }\right) $, we obtain the compact embedding 
\[
W_{r}\hookrightarrow L_{K}^{q_{1}}+L_{K}^{q_{2}}\qquad \text{for every }%
q_{1},q_{2}\text{\ satisfying (\ref{EX(SWW2):q1}) and (\ref{EX(SWW2):q2}).} 
\]
If furthermore $a,b,b_{0}$ are such that 
\begin{equation}
\frac{N+b}{N-p}<\frac{p(N-1)+pb_{0}-a}{p\left( N-1\right) +a(p-1)},
\label{EX1: further}
\end{equation}
then we can take $q_{1}=q_{2}$ and we get the compact embedding 
\begin{equation}
W_{r}\hookrightarrow L_{K}^{q}\qquad \text{for}\quad \max \left\{ 1,p\frac{%
N+b_{0}}{N+a},p\frac{N+b}{N-p}\right\} <q<p\frac{p(N-1)+pb_{0}-a}{p\left(
N-1\right) +a(p-1)}.  \label{EX1: Lq}
\end{equation}
Observe that the potentials $V$ and $K$ behave as a power and thus they fall
into the classes considered in \cite{Su-Wang-Will-p,SuTian12}. In
particular, the results of \cite{Su-Wang-Will-p} provide the compact
embedding 
\begin{equation}
W_{r}\hookrightarrow L_{K}^{q}\qquad \text{for}\quad \max \left\{ p,p\frac{%
N+b}{N-p}\right\} =:\underline{q}<q<\overline{q}:=p\frac{p(N-1)+pb_{0}-a}{%
p\left( N-1\right) +a(p-1)}.  \label{EX1: ssw}
\end{equation}
This requires condition (\ref{EX1: further}), which amounts to $\underline{q}%
<\overline{q}$, and no compact embedding is found in \cite{Su-Wang-Will-p}
if (\ref{EX1: further}) fails. Moreover, our result improves (\ref{EX1: ssw}%
) even if (\ref{EX1: further}) holds. Indeed, $b_{0}>a$ and $N+a<0$ imply $%
\frac{N+b_{0}}{N+a}<1$ and thus one has 
\[
\underline{q}=\left\{ 
\begin{array}{ll}
p\frac{N+b}{N-p}=\max \left\{ 1,p\frac{N+b_{0}}{N+a},p\frac{N+b}{N-p}%
\right\} \smallskip & \text{if }b\geq -p \\ 
p>\max \left\{ 1,p\frac{N+b_{0}}{N+a},p\frac{N+b}{N-p}\right\} & \text{if }%
b<-p
\end{array}
\right. , 
\]
so that (\ref{EX1: Lq}) is exactly (\ref{EX1: ssw}) if $b\geq -p$ and it is
better if $b<-p$. This last case actually concerns exponents less than $p$,
so it should be also compared with the results of \cite{SuTian12}, where,
setting 
\[
b_{1}:=\frac{p\left( N-1\right) +a(p-1)}{p^{2}}-N,\quad b_{2}:=\frac{p\left(
N-1\right) +a(p-1)}{p}-N,\quad b_{3}:=\frac{N-p}{p}-N 
\]
(notice that $-N<b_{1}<b_{2}<b_{3}<-p$ for $a$ as in (\ref{EX1: V})) and 
\begin{equation}
\underline{q}^{\prime }:=\left\{ 
\begin{array}{ll}
p\frac{N+b}{N-p}\smallskip & \text{if }b\in \left[ b_{3},-p\right) \\ 
\frac{p^{2}(N+b)}{p\left( N-1\right) +a(p-1)}\quad & \text{if }b\in \left[
b_{1},b_{2}\right)
\end{array}
\right. ,\qquad \overline{q}^{\prime }:=\left\{ 
\begin{array}{ll}
p\frac{N+b_{0}}{N-p} & \text{if }b_{0}\in \left( b_{3},-p\right] \\ 
\frac{p^{2}(N+b_{0})}{p\left( N-1\right) +a(p-1)}\quad & \text{if }b_{0}\in
\left( b_{1},b_{2}\right]
\end{array}
\right. ,  \label{EX1: st-qs}
\end{equation}
the authors find the compact embedding 
\begin{equation}
W_{r}\hookrightarrow L_{K}^{q}\qquad \text{for}\quad \underline{q}^{\prime
}<q<\overline{q}^{\prime }.  \label{EX1: st}
\end{equation}
Our result (\ref{EX1: Lq})-(\ref{EX1: further}) extends (\ref{EX1: st}) in
three directions. First, (\ref{EX1: st}) requires that $\underline{q}%
^{\prime }$ and $\overline{q}^{\prime }$ are defined, i.e., $b$ and $b_{0}$
lie in the intervals considered in (\ref{EX1: st-qs}), while (\ref{EX1: Lq})
and (\ref{EX1: further}) do not need such a restriction, also covering cases
of $b\in \left( -\infty ,b_{1}\right) \cup \left[ b_{2},b_{3}\right) $ or $%
b_{0}\in \left( a,b_{1}\right] \cup \left( b_{2},b_{3}\right] $ (take for
instance $b_{0}>a$ arbitrary and $b$ small enough to satisfy (\ref{EX1:
further})). Moreover, (\ref{EX1: st}) asks for the further condition $%
\underline{q}^{\prime }<\overline{q}^{\prime }$, which can be false even if $%
\underline{q}^{\prime }$ and $\overline{q}^{\prime }$ are defined (take for
instance $b=b_{0}\in \left( b_{3},-p\right) $ or $b=b_{0}\in \left(
b_{1},b_{2}\right) $, which give $\underline{q}^{\prime }=\overline{q}%
^{\prime }$), while condition (\ref{EX1: further}) does not. Actually, as
soon as $\underline{q}^{\prime }$ and $\overline{q}^{\prime }$ are defined,
one has $b<-p$ and $b_{0}>-N$, which imply 
\[
\frac{N+b}{N-p}-\frac{p(N-1)+pb_{0}-a}{p\left( N-1\right) +a(p-1)}<1+\frac{%
p+a}{p\left( N-1\right) +a(p-1)}=p\frac{N+a}{p\left( N-1\right) +a(p-1)}<0. 
\]
Finally, setting for brevity 
\[
\underline{q}^{\prime \prime }:=\max \left\{ 1,p\frac{N+b_{0}}{N+a},p\frac{%
N+b}{N-p}\right\} , 
\]
some computations (which we leave to the reader) show that, whenever $%
\underline{q}^{\prime }$ and $\overline{q}^{\prime }$ are defined, one has 
\[
\underline{q}^{\prime }=\left\{ 
\begin{array}{ll}
p\frac{N+b}{N-p}=\underline{q}^{\prime \prime }\medskip & \text{if }b\in
\left[ b_{3},-p\right) \\ 
\frac{p^{2}(N+b)}{p\left( N-1\right) +a(p-1)}=1=\underline{q}^{\prime \prime
}~~\medskip & \text{if }b=b_{1} \\ 
\frac{p^{2}(N+b)}{p\left( N-1\right) +a(p-1)}>1=\underline{q}^{\prime \prime
} & \text{if }b\in \left( b_{1},b_{2}\right)
\end{array}
\right. \quad \text{and}\quad \overline{q}^{\prime }<p\frac{p(N-1)+pb_{0}-a}{%
p\left( N-1\right) +a(p-1)}. 
\]
This shows that (\ref{EX1: Lq}) always gives a wider range of exponents $q$
than (\ref{EX1: st}).
\end{example}

\section{Existence and multiplicity results \label{SEC: ex}}

Let $1<p<N$. In this section we state our existence and multiplicity results
about radial weak solutions to the equation 
\begin{equation}
-\triangle _{p}u+V\left( \left| x\right| \right) |u|^{p-1}u=g\left( \left|
x\right| ,u\right) \quad \text{in }\mathbb{R}^{N},  \label{EQg}
\end{equation}
i.e., functions $u\in W_{r}$ such that 
\begin{equation}
\int_{\mathbb{R}^{N}}|\nabla u|^{p-2}\nabla u\cdot \nabla h\,dx+\int_{%
\mathbb{R}^{N}}V\left( \left| x\right| \right) |u|^{p-2}uh\,dx=\int_{%
\mathbb{R}^{N}}g\left( \left| x\right| ,u\right) h\,dx\quad \text{for all }%
h\in W,  \label{weak solution}
\end{equation}
where $V$ is a potential satisfying $\left( \mathbf{V}\right) $ and $W$ and $%
W_{r}$ are the Banach spaces defined in (\ref{spaces}), equipped with the
uniformly convex standard norm given by (\ref{norm}). As concerns the
nonlinearity, we assume that $g:\mathbb{R}_{+}\times \mathbb{R}\rightarrow \mathbb{R}$
is a Carath\'{e}odory function such that

\begin{itemize}
\item[$\left( g_{0}\right) $]  the linear operator $h\mapsto \int_{\mathbb{R}%
^{N}}g\left( \left| x\right| ,0\right) h\,dx$ is continuous on $W_{r}$
\end{itemize}

\noindent and that there exist $f\in C\left( \mathbb{R};\mathbb{R}\right) $ and a
function $K$ satisfying $\left( \mathbf{K}\right) $ such that:

\begin{itemize}
\item[$\left( g\right) $]  $\left| g\left( r,t\right) -g\left( r,0\right)
\right| \leq K\left( r\right) \left| f\left( t\right) \right| $ for almost
every $r>0$ and all $t\in \mathbb{R}.$
\end{itemize}

\noindent The model cases of $g$ we have in mind are, of course, $g\left(
r,t\right) =K\left( r\right) f\left( t\right) $ or $g\left( r,t\right)
=K\left( r\right) f\left( t\right) +Q\left( r\right) $ with $Q=g\left( \cdot
,0\right) $ such that $\left( g_{0}\right) $ holds (see Remark \ref
{RMK:operator}). On the function $f$ we will also require the following
condition (see also Remarks \ref{RMK:thm:ex}.\ref{RMK:thm:ex-2} and \ref
{RMK:ex-sub}.\ref{RMK:ex-sub-2}), where $q_{1},q_{2}$ will be specified each
time:

\begin{itemize}
\item[$\left( f_{q_{1},q_{2}}\right) $]  $\exists M>0$ such that $\left|
f\left( t\right) \right| \leq M\min \left\{ \left| t\right|
^{q_{1}-1},\left| t\right| ^{q_{2}-1}\right\} $ for all $t\in \mathbb{R}.$
\end{itemize}

\noindent Observe that, if $q_{1}\neq q_{2}$, the double-power growth
condition $\left( f_{q_{1},q_{2}}\right) $ is more stringent than the more
usual single-power one, since it implies $\sup_{t>0}\,\left| f\left(
t\right) \right| /t^{q-1}<+\infty $ for $q=q_{1}$, $q=q_{2}$ and every $q$
in between. On the other hand, we will never require $q_{1}\neq q_{2}$ in $%
\left( f_{q_{1},q_{2}}\right) $, so that our results will also concern
single-power nonlinearities as long as we can take $q_{1}=q_{2}$ (see
Example \ref{EX:eq} below).

\begin{rem}
\label{RMK:operator}Of course assumption $\left( g_{0}\right) $ will be
relevant only if $g\left( \cdot ,0\right) \neq 0$ (meaning that $g\left(
\cdot ,0\right) $ does not vanish almost everywhere). In this case, the
radial estimates satisfied by the $W_{r}$ mappings (see Lemmas 4.2 and 4.3
of \cite{BGR-p.arxiv}) provide simple explicit conditions ensuring
assumption $\left( g_{0}\right) $, which turns out to be fulfilled if $%
g\left( \left| \cdot \right| ,0\right) $ belongs to $L_{\mathrm{loc}}^{1}(%
\mathbb{R}^{N}\setminus \left\{ 0\right\} )$ and satisfies suitable decay (or
growth) conditions at zero and at infinity. On the other hand, it is not
difficult to find explicit conditions on $g\left( \cdot ,0\right) $ ensuring 
$\left( g_{0}\right) $ even on the whole space $W$, for example $g\left(
\cdot ,0\right) \in L^{p/(p-1)}(\mathbb{R}_{+},r^{N+1/(p-1)}dr)$. Indeed, this
means $g\left( \left| \cdot \right| ,0\right) \in L^{p/(p-1)}(\mathbb{R}%
^{N},\left| x\right| ^{p/(p-1)}dx)$ and thus it implies 
\[
\left| \int_{\mathbb{R}^{N}}g\left( \left| x\right| ,0\right) h\,dx\right| \leq
\left( \int_{\mathbb{R}^{N}}\left| g\left( \left| x\right| ,0\right) \right| ^{%
\frac{p}{p-1}}\left| x\right| ^{\frac{p}{p-1}}dx\right) ^{\frac{p-1}{p}%
}\left( \int_{\mathbb{R}^{N}}\frac{\left| h\right| ^{p}}{\left| x\right| ^{p}}%
dx\right) ^{\frac{1}{p}}\leq \left( \mathrm{const.}\right) \left\| h\right\| 
\]
for all $h\in W\hookrightarrow D^{1,p}(\mathbb{R}^{N})$, by H\"{o}lder and
Hardy inequalities. Other conditions ensuring the same result are $g\left(
\cdot ,0\right) \in L^{pN/\left( pN-N+p\right) }(\mathbb{R}_{+},r^{N-1}dr)$ or $%
V^{-1/p}g\left( \cdot ,0\right) \in L^{p/(p-1)}(\mathbb{R}_{+},r^{N-1}dr)$.
\end{rem}

Set $G\left( r,t\right) :=\int_{0}^{t}g\left( r,s\right) ds$ and 
\begin{equation}
I\left( u\right) :=\frac{1}{p}\left\| u\right\| ^{p}-\int_{\mathbb{R}%
^{N}}G\left( \left| x\right| ,u\right) dx.  \label{I:=}
\end{equation}
From the continuous embedding result of Theorem \ref{THM(cpt)} and the
results of \cite{BPR} about Nemytski\u{\i} operators on the sum of Lebesgue
spaces, we have that (\ref{I:=}) defines a $C^{1}$ functional on $W_{r}$
provided that there exist $q_{1},q_{2}>1$ such that $\left(
f_{q_{1},q_{2}}\right) $ and $\left( \mathcal{S}_{q_{1},q_{2}}^{\prime
}\right) $ hold. In this case, the Fr\'{e}chet derivative of $I$ at any $%
u\in W_{r}$ is given by 
\begin{equation}
I^{\prime }\left( u\right) h=\int_{\mathbb{R}^{N}}\left( |\nabla u|^{p-2}\nabla
u\cdot \nabla h+V\left( \left| x\right| \right) |u|^{p-2}uh\right) dx-\int_{%
\mathbb{R}^{N}}g\left( \left| x\right| ,u\right) h\,dx,\quad \forall h\in W_{r}
\label{PROP:diff: I'(u)h=}
\end{equation}
and therefore the critical points of $I:W_{r}\rightarrow \mathbb{R}$ satisfy (%
\ref{weak solution}) for all $h\in W_{r}$. Our first result shows that such
critical points are actually weak solutions to equation (\ref{EQg}),
provided that the following slightly stronger version of condition $\left( 
\mathcal{S}_{q_{1},q_{2}}^{\prime }\right) $ holds: 
\begin{equation}
\mathcal{R}_{0}\left( q_{1},R_{1}\right) <\infty \quad \text{and}\quad 
\mathcal{R}_{\infty }\left( q_{2},R_{2}\right) <\infty \quad \text{for some }%
R_{1},R_{2}>0.  
\tag*{$\left( {\cal R}_{q_{1},q_{2}}^{\prime }\right) $}
\end{equation}
Observe that the classical Palais' Principle of Symmetric Criticality \cite
{Palais} does not apply in this case, because we do not know whether or not $%
I$ is differentiable, not even well defined, on the whole space $W$.

\begin{prop}
\label{PROP:symm-crit}Assume $s>\frac{Np}{N(p-1)+p}$ in condition $(\mathbf{K%
})$ and assume that $\left( g_{0}\right) $ holds on the whole space $W$ (cf.
Remark \ref{RMK:operator}). Assume furthermore that there exist $%
q_{1},q_{2}>1$ such that $\left( f_{q_{1},q_{2}}\right) $ and $\left( 
\mathcal{R}_{q_{1},q_{2}}^{\prime }\right) $ hold. Then every critical point
of $I:W_{r}\rightarrow \mathbb{R}$ is a weak solution to equation (\ref{EQg})).
\end{prop}

By Proposition \ref{PROP:symm-crit}, the problem of radial weak solutions to
(\ref{EQg}) reduces to the study of the critical points of $%
I:W_{r}\rightarrow \mathbb{R}$, which is the aim of our next results.

Concerning the case of super $p$-linear nonlinearities, we will prove the
following existence and multiplicity theorems.

\begin{thm}
\label{THM:ex}Assume $g\left( \cdot ,0\right) =0$ and assume that there
exist $q_{1},q_{2}>p$ such that $\left( f_{q_{1},q_{2}}\right) $ and $\left( 
\mathcal{S}_{q_{1},q_{2}}^{\prime \prime }\right) $ hold. Assume furthermore
that $g$ satisfies:

\begin{itemize}
\item[$\left( g_{1}\right) $]  $\exists \theta >p$ such that $0\leq \theta
G\left( r,t\right) \leq g\left( r,t\right) t$ for almost every $r>0$ and all 
$t\geq 0;$

\item[$\left( g_{2}\right) $]  $\exists t_{0}>0$ such that $G\left(
r,t_{0}\right) >0$ for almost every $r>0.$
\end{itemize}

\noindent If $K\left( \left| \cdot \right| \right) \in L^{1}(\mathbb{R}^{N})$,
we can replace assumptions $\left( g_{1}\right) $-$\left( g_{2}\right) $
with:

\begin{itemize}
\item[$\left( g_{3}\right) $]  $\exists \theta >p$ and $\exists t_{0}>0$
such that $0<\theta G\left( r,t\right) \leq g\left( r,t\right) t$ for almost
every $r>0$ and all $t\geq t_{0}.$
\end{itemize}

\noindent Then the functional $I:W_{r}\rightarrow \mathbb{R}$ has a nonnegative
critical point $u\neq 0$.
\end{thm}

\begin{rem}
\label{RMK:thm:ex}\quad

\begin{enumerate}
\item  \label{RMK:thm:ex-1}Assumptions $\left( g_{1}\right) $ and $\left(
g_{2}\right) $ imply $\left( g_{3}\right) $, so that, in Theorem \ref{THM:ex}%
, the information $K\left( \left| \cdot \right| \right) \in L^{1}(\mathbb{R}%
^{N})$ actually allows weaker hypotheses on the nonlinearity.

\item  \label{RMK:thm:ex-2}In Theorem \ref{THM:ex}, assumptions $\left(
g\right) $ and $\left( f_{q_{1},q_{2}}\right) $ need only to hold for $t\geq
0$. Indeed, all the hypotheses of the theorem still hold true if we replace $%
g\left( r,t\right) $ with $\chi _{\mathbb{R}_{+}}\left( t\right) g\left(
r,t\right) $ ($\chi _{\mathbb{R}_{+}}$ is the characteristic function of $\mathbb{R%
}_{+}$) and this can be done without restriction since the theorem concerns 
\emph{nonnegative} critical points.
\end{enumerate}
\end{rem}

\begin{thm}
\label{THM:mult}Assume that there exist $q_{1},q_{2}>p$ such that $\left(
f_{q_{1},q_{2}}\right) $ and $\left( \mathcal{S}_{q_{1},q_{2}}^{\prime
\prime }\right) $ hold. Assume furthermore that:

\begin{itemize}
\item[$\left( g_{4}\right) $]  $\exists m>0$ such that $G\left( r,t\right)
\geq mK\left( r\right) \min \left\{ t^{q_{1}},t^{q_{2}}\right\} $ for almost
every $r>0$ and all $t\geq 0;$

\item[$\left( g_{5}\right) $]  $g\left( r,t\right) =-g\left( r,-t\right) $
for almost every $r>0$ and all $t\geq 0.$
\end{itemize}

\noindent Finally, assume that $g$ satisfies $\left( g_{1}\right) $, or that 
$K\left( \left| \cdot \right| \right) \in L^{1}(\mathbb{R}^{N})$ and $g$
satisfies $\left( g_{3}\right) $. Then the functional $I:W_{r}\rightarrow 
\mathbb{R}$ has a sequence of critical points $\left\{ u_{n}\right\} $ such
that $I\left( u_{n}\right) \rightarrow +\infty $.
\end{thm}

\begin{rem}
The condition $g\left( \cdot ,0\right) =0$ is implicit in Theorem \ref
{THM:mult} (and in Theorem \ref{THM:mult-sub} below), as it follows from
assumption $\left( g_{5}\right) $.
\end{rem}

As to sub $p$-linear nonlinearities, we will prove the following results,
where we also consider the case $g\left( \cdot ,0\right) \neq 0$.

\begin{thm}
\label{THM:ex-sub}Assume that there exist $q_{1},q_{2}\in \left( 1,p\right) $
such that $\left( f_{q_{1},q_{2}}\right) $ and $\left( \mathcal{S}%
_{q_{1},q_{2}}^{\prime \prime }\right) $ hold. Assume furthermore that $g$
satisfies $\left( g_{0}\right) $ and at least one of the following
conditions:

\begin{itemize}
\item[$\left( g_{6}\right) $]  $\exists \theta <p$ and $\exists t_{0},m>0$
such that $G\left( r,t\right) \geq mK\left( r\right) t^{\theta }$ for almost
every $r>0$ and all $0\leq t\leq t_{0};$

\item[$\left( g_{7}\right) $]  $g\left( \cdot ,0\right) $ does not vanish
almost everywhere in $\left( r_{1},r_{2}\right) .$
\end{itemize}

\noindent If $\left( g_{7}\right) $ holds, we also allow the case $\max
\left\{ q_{1},q_{2}\right\} =p>\min \left\{ q_{1},q_{2}\right\} >1$. Then
there exists $u\neq 0$ such that 
\[
I\left( u\right) =\min_{v\in W_{r}}I\left( v\right) . 
\]
\end{thm}

If $g\left( \cdot ,t\right) \geq 0$ almost everywhere for all $t<0$, the
minimizer $u$ of Theorem \ref{THM:ex-sub} is nonnegative, since a standard
argument shows that all the critical points of $I$ are nonnegative (test $%
I^{\prime }\left( u\right) $ with the negative part $u_{-}$ and get $%
I^{\prime }\left( u\right) u_{-}=-\left\| u_{-}\right\| ^{p}=0$). The next
corollary gives a nonnegative critical point just requiring $g\left( \cdot
,0\right) \geq 0$.

\begin{cor}
\label{COR:ex-sub}Assume the same hypotheses of Theorem \ref{THM:ex-sub}. If 
$g\left( \cdot ,0\right) \geq 0$ almost everywhere, then $I:W_{r}\rightarrow 
\mathbb{R}$ has a nonnegative critical point $\widetilde{u}\neq 0$ satisfying 
\begin{equation}
I\left( \widetilde{u}\right) =\min_{u\in W_{r},\,u\geq 0}I\left( u\right) .
\label{COR:ex-sub: th}
\end{equation}
\end{cor}

\begin{rem}
\label{RMK:ex-sub}\quad

\begin{enumerate}
\item  In Theorem \ref{THM:ex-sub} and Corollary \ref{COR:ex-sub}, the case $%
\max \left\{ q_{1},q_{2}\right\} =p>\min \left\{ q_{1},q_{2}\right\} >1$
cannot be considered under assumption $\left( g_{6}\right) $, since $\left(
g_{6}\right) $ and $\left( f_{q_{1},q_{2}}\right) $ imply $\max \left\{
q_{1},q_{2}\right\} \leq \theta <p$.

\item  \label{RMK:ex-sub-2}Checking the proof, one sees that Corollary \ref
{COR:ex-sub} actually requires that assumptions $\left( g\right) $ and $%
\left( f_{q_{1},q_{2}}\right) $ hold only for $t\geq 0$, which is consistent
with the concern of the result about \emph{nonnegative} critical points.
\end{enumerate}
\end{rem}

\begin{thm}
\label{THM:mult-sub}Assume that there exist $q_{1},q_{2}\in \left(
1,p\right) $ such that $\left( f_{q_{1},q_{2}}\right) $ and $\left( \mathcal{%
S}_{q_{1},q_{2}}^{\prime \prime }\right) $ hold. Assume furthermore that $g$
satisfies $\left( g_{5}\right) $ and $\left( g_{6}\right) $. Then the
functional $I:W_{r}\rightarrow \mathbb{R}$ has a sequence of critical points $%
\left\{ u_{n}\right\} $ such that $I\left( u_{n}\right) <0$ and $I\left(
u_{n}\right) \rightarrow 0$.
\end{thm}

All the above existence and multiplicity results rely on assumption $\left( 
\mathcal{S}_{q_{1},q_{2}}^{\prime \prime }\right) $, which is quite abstract
but can be granted in concrete cases through Theorems \ref{THM0}, \ref{THM1}%
, \ref{THM2} and \ref{THM3}, which ensure the stronger condition $\left( 
\mathcal{R}_{q_{1},q_{2}}^{\prime \prime }\right) $ for suitable ranges of
exponents $q_{1}$ and $q_{2}$ by explicit conditions on the potentials. This
has been already discussed in Section \ref{SEC: comp} and exemplified in
Section \ref{SEC: examples}. Moreover, explicit conditions on a forcing term 
$g\left( \cdot ,0\right) \neq 0$ in order that $\left( g_{0}\right) $ holds
has been examined in Remark \ref{RMK:operator}, so here we limit ourselves
to give some basic examples of nonlinearities of the form $g\left(
r,t\right) =K\left( r\right) f\left( t\right) $ satisfying the assumptions
of our results and to apply them to a sample equation.

\begin{example}
\label{EX:f}Let $g\left( r,t\right) =K\left( r\right) f\left( t\right) $
with $K$ satisfying $\left( \mathbf{K}\right) $. The simplest $f\in C\left( 
\mathbb{R};\mathbb{R}\right) $ such that $\left( f_{q_{1},q_{2}}\right) $ holds is 
$f\left( t\right) =\min \left\{ \left| t\right| ^{q_{1}-2}t,\left| t\right|
^{q_{2}-2}t\right\} $, which also ensures $\left( g_{1}\right) $ if $%
q_{1},q_{2}>p$ (with $\theta =\min \left\{ q_{1},q_{2}\right\} $), and $%
\left( g_{6}\right) $ if $q_{1},q_{2}<p$ (with $\theta =\max \left\{
q_{1},q_{2}\right\} $). Another model example is 
\[
f\left( t\right) =\frac{\left| t\right| ^{q_{2}-2}t}{1+\left| t\right|
^{q_{2}-q_{1}}}\quad \text{with }1<q_{1}\leq q_{2},
\]
which ensures $\left( g_{1}\right) $ if $q_{1}>p$ (with $\theta =q_{1}$) and 
$\left( g_{6}\right) $ if $q_{2}<p$ (with $\theta =q_{2}$). Note that, in
both these cases, also $\left( g_{2}\right) $, $\left( g_{4}\right) $ and $%
\left( g_{5}\right) $ hold true. Moreover, both of these functions $f$
become $f\left( t\right) =\left| t\right| ^{q-2}t$ if $q_{1}=q_{2}=q$. Other
examples of nonlinearities $f$ ensuring $\left( f_{q_{1},q_{2}}\right) $ are 
\[
f\left( t\right) =\frac{\left| t\right| ^{q_{1}+q-1}-\left| t\right|
^{q_{2}-1}}{1+\left| t\right| ^{q}},\quad f\left( t\right) =\frac{\left|
t\right| ^{q_{2}-1+\varepsilon }}{1+\left| t\right|
^{q_{2}-q_{1}+2\varepsilon }}\ln \left| t\right| 
\]
(the latter extended at $0$ by continuity) with $1<q_{1}\leq q_{2}<q_{1}+q$
and $\varepsilon >0$, for which $\left( g_{1}\right) $ and $\left(
g_{6}\right) $ do not hold, but $\left( g_{3}\right) $ is granted if $q_{1}>p
$ and $\varepsilon $ is small enough.
\end{example}

\begin{example}
\label{EX:eq}Let $1<p<N$ and $a>0$, and consider the equation 
\begin{equation}
-\triangle _{p}u+\frac{e^{-a\left| x\right| }}{\left| x\right| ^{N}}%
|u|^{p-1}u=K\left( \left| x\right| \right) \left| u\right| ^{q-1}u\quad 
\text{in }\mathbb{R}^{N}  \label{power-eq}
\end{equation}
where $K:\mathbb{R}_{+}\rightarrow \mathbb{R}_{+}$ is a continuous function such
that $K\left( r\right) =O(r^{b_{0}})_{r\rightarrow 0^{+}}$ and $K\left(
r\right) =O(r^{b})_{r\rightarrow +\infty }$ with $b_{0}>-N$ and $b\in %
\mathbb{R}$. By Theorems \ref{THM1} and \ref{THM3} (applied with $\gamma
_{0}=N$, $\alpha _{0}=b_{0}$, $\alpha _{\infty }=b$, $\beta _{0}=\beta
_{\infty }=0$), condition $\left( \mathcal{S}_{q_{1},q_{2}}^{\prime \prime
}\right) $ hold if 
\begin{equation}
1<q_{1}<\overline{q}:=p\left( 1+\frac{N+b_{0}}{N-p}p\right) \quad \text{and}%
\quad q_{2}>\underline{q}:=\max \left\{ 1,p\frac{N+b}{N-p}\right\} .
\label{range}
\end{equation}
Note that $\overline{q}>p$, since $b_{0}>-N$. Then, by Theorem \ref{THM:ex},
Corollary \ref{COR:ex-sub} and Proposition \ref{PROP:symm-crit}, the
equation has a nonnegative radial weak solution in the following cases:

\begin{itemize}
\item  $b<p\left( N+b_{0}-1\right) $ and $\max \left\{ p,p\frac{N+b}{N-p}%
\right\} <q<\overline{q};$

\item  $b<-p$ and $\underline{q}<q<p.$
\end{itemize}

\noindent If $b\geq p\left( N+b_{0}-1\right) $, instead, we cannot pick $%
q_{1}=q_{2}>p$ in (\ref{range}) (since $\overline{q}\leq \underline{q}$),
nor $q_{1}=q_{2}<p$ (since $\underline{q}>p$). In this case our results do
not apply to equation (\ref{power-eq}), but they apply to the equation
\[
-\triangle _{p}u+\frac{e^{-a\left| x\right| }}{\left| x\right| ^{N}}%
|u|^{p-1}u=K\left( \left| x\right| \right) f\left( u\right) \quad \text{in }%
\mathbb{R}^{N}
\]
where $f$ is any nonlinearity ensuring $\left( g_{1}\right) $ and $\left(
g_{2}\right) $ and satisfying $\left( f_{q_{1},q_{2}}\right) $ (for instance
one of the first two functions considered in Example \ref{EX:f}), for which
Theorem \ref{THM:ex} and Proposition \ref{PROP:symm-crit} provide a
nonnegative radial weak solution if $p<q_{1}<\overline{q}$ and $q_{2}>%
\underline{q}$.
\end{example}

We end this section by observing that from the above results one can also
derive existence and multiplicity results for equation (\ref{EQ}) with
Dirichlet boundary conditions in bounded balls or exterior radial domains,
where a single-power growth condition on the nonlinearity is sufficient and,
respectively, only assumptions on $V$ and $K$ near the origin or at infinity
are needed. This can be done by suitably modifying the potentials $V$ and $K$%
, in order to reduce the Dirichlet problem to the problem in $\mathbb{R}^{N}$
(see \cite[Section 5]{BGR_II}). We leave the details to the interested
reader, as well as the precise statements of the results.

\section{Proofs \label{SEC:pf-general}}

This section is devoted to proof of the results of Section \ref{SEC: ex}, so
we keep the notation and assumptions of that section. We begin by recalling
the following lemma from \cite{BGR-p.arxiv}.

\begin{lem}[{\cite[Lemma 3.1]{BGR-p.arxiv}}]
\label{LEM:corona} Let $R>r>0$ and $1<q<\infty $. Then there exist $\tilde{C}%
=\tilde{C}\left( N,p,r,R,q,s\right) >0$ and $l=l\left( p,q,s\right) >0$ such
that $\forall u\in W_{r}$ one has 
\[
\int_{B_{R}\setminus B_{r}}K\left( \left| x\right| \right) \left| u\right|
^{q}dx\leq \tilde{C}\left\| K\left( \left| \cdot \right| \right) \right\|
_{L^{s}(B_{R}\setminus B_{r})}\left\| u\right\| ^{q-lp}\left(
\int_{B_{R}\setminus B_{r}}\left| u\right| ^{p}dx\right) ^{l}.
\]
Moreover, if $s>\frac{Np}{N(p-1)+p}$ in assumption $(\mathbf{K})$, then
there exists $\tilde{C}_{1}=\tilde{C}_{1}\left( N,p,r,R,q,s\right) >0$ such
that $\forall u\in W_{r}$ and $\forall h\in W$ one has 
\[
\frac{\int_{B_{R}\setminus B_{r}}K\left( \left| x\right| \right) \left|
u\right| ^{q-1}\left| h\right| dx}{\tilde{C}_{1}\left\| K\left( \left| \cdot
\right| \right) \right\| _{L^{s}(B_{R}\setminus B_{r})}}\leq \left\{ 
\begin{array}{ll}
\left( \int_{B_{R}\setminus B_{r}}\left| u\right| ^{p}dx\right) ^{\frac{q-1}{%
p}}\left\| h\right\| \medskip  & \text{if }q\leq \tilde{q} \\ 
\left( \int_{B_{R}\setminus B_{r}}\left| u\right| ^{p}dx\right) ^{\frac{%
\tilde{q}-1}{p}}\left\| u\right\| ^{q-\tilde{q}}\left\| h\right\| \quad
\medskip  & \text{if }q>\tilde{q}
\end{array}
\right. 
\]
where $\tilde{q}:=p\left( 1+\frac{1}{N}-\frac{1}{s}\right) $ (note that $s>%
\frac{Np}{N(p-1)+p}$ implies $\tilde{q}>1$).
\end{lem}

\noindent \textbf{Proof of Proposition \ref{PROP:symm-crit}.}\quad Let $u\in
W_{r}$. By the monotonicity of $\mathcal{R}_{0}$ and $\mathcal{R}_{\infty }$%
, it is not restrictive to assume $R_{1}<R_{2}$ in hypothesis $\left( 
\mathcal{R}_{q_{1},q_{2}}\right) $. So, by Lemma \ref{LEM:corona}, there
exists a constant $C>0$ (dependent on $u$) such that for all $h\in W$ we
have 
\[
\int_{B_{R_{2}}\setminus B_{R_{1}}}K\left( \left| x\right| \right) \left|
u\right| ^{q_{1}-1}\left| h\right| dx\leq C\left\| h\right\| 
\]
and therefore, by $\left( g\right) $ and $\left( f_{q_{1},q_{2}}\right) $, 
\begin{eqnarray*}
\int_{\mathbb{R}^{N}}\left| \,g\left( \left| x\right| ,u\right) -g\left( \left|
x\right| ,0\right) \,\right| \left| h\right| dx &\leq &M\int_{\mathbb{R}%
^{N}}K\left( \left| x\right| \right) \min \{\left| u\right|
^{q_{1}-1},\left| u\right| ^{q_{2}-1}\}\left| h\right| dx \\
&\leq &M\left( \int_{B_{R_{1}}}K\left( \left| x\right| \right) \left|
u\right| ^{q_{1}-1}\left| h\right| dx+\int_{B_{R_{2}}^{c}}K\left( \left|
x\right| \right) \left| u\right| ^{q_{2}-1}\left| h\right| dx+C\left\|
h\right\| \right) \\
&\leq &M\left( \left\| u\right\| ^{q_{1}-1}\mathcal{R}_{0}\left(
q_{1},R_{1}\right) +\left\| u\right\| ^{q_{2}-1}\mathcal{R}_{\infty }\left(
q_{2},R_{2}\right) +C\right) \left\| h\right\| .
\end{eqnarray*}
Together with the assumption on the continuity on $W$ of the operator $%
h\mapsto \int_{\mathbb{R}^{N}}g\left( \left| x\right| ,0\right) h\,dx$, this
gives that the linear operator 
\[
T\left( u\right) h:=\int_{\mathbb{R}^{N}}\left( |\nabla u|^{p-2}\nabla u\cdot
\nabla h+V\left( \left| x\right| \right) |u|^{p-2}uh\right) dx-\int_{\mathbb{R}%
^{N}}g\left( \left| x\right| ,u\right) h\,dx 
\]
is well defined and continuous on $W$. Hence, by uniform convexity, there
exists a unique $\tilde{u}\in W$ such that $T\left( u\right) \tilde{u}%
=\left\| \tilde{u}\right\| ^{2}=\left\| T\left( u\right) \right\|
_{W^{\prime }}^{2}$. Denoting by $O\left( N\right) $ the orthogonal group of 
$\mathbb{R}^{N}$, by means of obvious changes of variables it is easy to see
that for every $h\in W$ one has 
\[
T\left( u\right) h\left( S\cdot \right) =T\left( u\right) h\quad \text{and}%
\quad \left\| h\left( S\cdot \right) \right\| =\left\| h\right\| \quad \text{%
for all }S\in O\left( N\right) , 
\]
whence, applying with $h=\tilde{u}$, one deduces $\tilde{u}\left( S\cdot
\right) =\tilde{u}$ by uniqueness. This means $\tilde{u}\in W_{r}$, so that,
if $T\left( u\right) h=0$ for all $h\in W_{r}$, one gets $T\left( u\right) 
\tilde{u}=0$ and hence $\left\| T\left( u\right) \right\| _{W^{\prime }}=0$.%
\endproof%
\bigskip

For future reference, we observe here that, by assumption $\left( g\right) $%
, if $\left( f_{q_{1},q_{2}}\right) $ holds then there exists $\tilde{M}>0$
such that for almost every $r>0$ and all $t\in \mathbb{R}$ one has 
\begin{equation}
\left| G\left( r,t\right) -g\left( r,0\right) t\right| \leq \tilde{M}K\left(
r\right) \min \left\{ \left| t\right| ^{q_{1}},\left| t\right|
^{q_{2}}\right\} .  \label{G_pq}
\end{equation}

\begin{lem}
\label{LEM:MPgeom}Assume $\left( g_{0}\right) $ and let $L_{0}$ be the norm
of the operator therein. If there exist $q_{1},q_{2}>1$ such that $\left(
f_{q_{1},q_{2}}\right) $ and $\left( \mathcal{S}_{q_{1},q_{2}}^{\prime
}\right) $ hold, then there exist two constants $c_{1},c_{2}>0$ such that 
\begin{equation}
I\left( u\right) \geq \frac{1}{p}\left\| u\right\| ^{p}-c_{1}\left\|
u\right\| ^{q_{1}}-c_{2}\left\| u\right\| ^{q_{2}}-L_{0}\left\| u\right\|
\qquad \text{for all }u\in W_{r}.  \label{LEM:MPgeom: th}
\end{equation}
If $\left( \mathcal{S}_{q_{1},q_{2}}^{\prime \prime }\right) $ also holds,
then $\forall \varepsilon >0$ there exist two constants $c_{1}\left(
\varepsilon \right) ,c_{2}\left( \varepsilon \right) >0$ such that (\ref
{LEM:MPgeom: th}) holds both with $c_{1}=\varepsilon $, $c_{2}=c_{2}\left(
\varepsilon \right) $ and with $c_{1}=c_{1}\left( \varepsilon \right) $, $%
c_{2}=\varepsilon $.
\end{lem}

\proof%
Let $i\in \left\{ 1,2\right\} $. By the monotonicity of $\mathcal{S}_{0}$
and $\mathcal{S}_{\infty }$, it is not restrictive to assume $R_{1}<R_{2}$
in hypothesis $\left( \mathcal{S}_{q_{1},q_{2}}^{\prime }\right) $. Then, by
Lemma \ref{LEM:corona} and the continuous embedding $W\hookrightarrow L_{%
\mathrm{loc}}^{p}(\mathbb{R}^{N})$, there exists a constant $%
c_{R_{1},R_{2}}^{\left( i\right) }>0$ such that for all $u\in W_{r}$ we have 
\[
\int_{B_{R_{2}}\setminus B_{R_{1}}}K\left( \left| x\right| \right) \left|
u\right| ^{q_{i}}dx\leq c_{R_{1},R_{2}}^{\left( i\right) }\left\| u\right\|
^{q_{i}}. 
\]
Therefore, by (\ref{G_pq}) and the definitions of $\mathcal{S}_{0}$ and $%
\mathcal{S}_{\infty }$, we obtain 
\begin{eqnarray}
&&\left| \int_{\mathbb{R}^{N}}G\left( \left| x\right| ,u\right) dx\right| \leq
\int_{\mathbb{R}^{N}}\left| G\left( \left| x\right| ,u\right) -g\left( \left|
x\right| ,0\right) u\right| dx+\left| \int_{\mathbb{R}^{N}}g\left( \left|
x\right| ,0\right) u\,dx\right|  \nonumber \\
&\leq &\tilde{M}\int_{\mathbb{R}^{N}}K\left( \left| x\right| \right) \min
\left\{ \left| u\right| ^{q_{1}},\left| u\right| ^{q_{2}}\right\}
dx+L_{0}\left\| u\right\|  \nonumber \\
&\leq &\tilde{M}\left( \int_{B_{R_{1}}}K\left( \left| x\right| \right)
\left| u\right| ^{q_{1}}dx+\int_{B_{R_{2}}^{c}}K\left( \left| x\right|
\right) \left| u\right| ^{q_{2}}dx+\int_{B_{R_{2}}\setminus
B_{R_{1}}}K\left( \left| x\right| \right) \left| u\right| ^{q_{i}}dx\right)
+L_{0}\left\| u\right\|  \nonumber \\
&\leq &\tilde{M}\left( \left\| u\right\| ^{q_{1}}\mathcal{S}_{0}\left(
q_{1},R_{1}\right) +\left\| u\right\| ^{q_{2}}\mathcal{S}_{\infty }\left(
q_{2},R_{2}\right) +c_{R_{1},R_{2}}^{\left( i\right) }\left\| u\right\|
^{q_{i}}\right) +L_{0}\left\| u\right\|  \label{LEM:MPgeom: pf} \\
&=&c_{1}\left\| u\right\| ^{q_{1}}+c_{2}\left\| u\right\|
^{q_{2}}+L_{0}\left\| u\right\| ,  \nonumber
\end{eqnarray}
with obvious definition of the constants $c_{1}$ and $c_{2}$, independent of 
$u$. This yields (\ref{LEM:MPgeom: th}). If $\left( \mathcal{S}%
_{q_{1},q_{2}}^{\prime \prime }\right) $ also holds, then $\forall
\varepsilon >0$ we can fix $R_{1,\varepsilon }<R_{2,\varepsilon }$ such that 
$\tilde{M}\mathcal{S}_{0}\left( q_{1},R_{1,\varepsilon }\right) <\varepsilon 
$ and $\tilde{M}\mathcal{S}_{\infty }\left( q_{2},R_{2,\varepsilon }\right)
<\varepsilon $, so that inequality (\ref{LEM:MPgeom: pf}) becomes 
\[
\left| \int_{\mathbb{R}^{N}}G\left( \left| x\right| ,u\right) dx\right| \leq
\varepsilon \left\| u\right\| ^{q_{1}}+\varepsilon \left\| u\right\|
^{q_{2}}+c_{R_{1,\varepsilon },R_{2,\varepsilon }}^{\left( i\right) }\left\|
u\right\| ^{q_{i}}+L_{0}\left\| u\right\| . 
\]
The result then ensues by taking $i=2$ and $c_{2}\left( \varepsilon \right)
=\varepsilon +c_{R_{1,\varepsilon },R_{2,\varepsilon }}^{\left( 2\right) }$,
or $i=1$ and $c_{1}\left( \varepsilon \right) =\varepsilon
+c_{R_{1,\varepsilon },R_{2,\varepsilon }}^{\left( 1\right) }$.%
\endproof%
\bigskip

Henceforth, we will assume that the hypotheses of Theorem \ref{THM:ex} also
include the following condition: 
\begin{equation}
g\left( r,t\right) =0\quad \text{for all }r>0~\text{and~}t<0.  \label{g=0}
\end{equation}
This can be done without restriction, since the theorem concerns \emph{%
nonnegative} critical points and all its assumptions still hold true if we
replace $g\left( r,t\right) $ with $g\left( r,t\right) \chi _{\mathbb{R}%
_{+}}\left( t\right) $ ($\chi _{\mathbb{R}_{+}}$ is the characteristic function
of $\mathbb{R}_{+}$).

\begin{lem}
\label{LEM:PS}Under the assumptions of each of Theorems \ref{THM:ex}
(including (\ref{g=0})) and \ref{THM:mult}, the functional $%
I:W_{r}\rightarrow \mathbb{R}$ satisfies the Palais-Smale condition.
\end{lem}

\proof%
By (\ref{g=0}) and $\left( g_{5}\right) $ respectively, under the
assumptions of each of Theorems \ref{THM:ex} and \ref{THM:mult} we have that
either $g$ satisfies $\left( g_{1}\right) $ for all $t\in \mathbb{R}$, or $%
K\left( \left| \cdot \right| \right) \in L^{1}(\mathbb{R}^{N})$ and $g$
satisfies 
\begin{equation}
\theta G\left( r,t\right) \leq g\left( r,t\right) t\quad \text{for almost
every }r>0~\text{and all~}\left| t\right| \geq t_{0}.  \label{LEM:PS: AR}
\end{equation}
Let $\left\{ u_{n}\right\} $ be a sequence in $W_{r}$ such that $\left\{
I\left( u_{n}\right) \right\} $ is bounded and $I^{\prime }\left(
u_{n}\right) \rightarrow 0$ in $W_{r}^{\prime }$. Hence 
\[
\frac{1}{p}\left\| u_{n}\right\| ^{p}-\int_{\mathbb{R}^{N}}G\left( \left|
x\right| ,u_{n}\right) dx=O\left( 1\right) \quad \text{and}\quad \left\|
u_{n}\right\| ^{p}-\int_{\mathbb{R}^{N}}g\left( \left| x\right| ,u_{n}\right)
u_{n}dx=o\left( 1\right) \left\| u_{n}\right\| . 
\]
If $g$ satisfies $\left( g_{1}\right) $, then we get 
\[
\frac{1}{p}\left\| u_{n}\right\| ^{p}+O\left( 1\right) =\int_{\mathbb{R}%
^{N}}G\left( \left| x\right| ,u_{n}\right) dx\leq \frac{1}{\theta }\int_{%
\mathbb{R}^{N}}g\left( \left| x\right| ,u_{n}\right) u_{n}dx=\frac{1}{\theta }%
\left\| u_{n}\right\| ^{p}+o\left( 1\right) \left\| u_{n}\right\| , 
\]
which implies that $\left\{ \left\| u_{n}\right\| \right\} $ is bounded
since $\theta >p$. If $K\left( \left| \cdot \right| \right) \in L^{1}(\mathbb{R}%
^{N})$ and $g$ satisfies (\ref{LEM:PS: AR}), then we slightly modify the
argument: we have 
\[
\int_{\left\{ \left| u_{n}\right| \geq t_{0}\right\} }g\left( \left|
x\right| ,u_{n}\right) u_{n}dx\leq \int_{\mathbb{R}^{N}}g\left( \left| x\right|
,u_{n}\right) u_{n}dx+\int_{\left\{ \left| u_{n}\right| <t_{0}\right\}
}\left| g\left( \left| x\right| ,u_{n}\right) u_{n}\right| dx 
\]
where (thanks to $\left( g\right) $ and $\left( f_{q_{1},q_{2}}\right) $) 
\[
\int_{\left\{ \left| u_{n}\right| <t_{0}\right\} }\left| g\left( \left|
x\right| ,u_{n}\right) u_{n}\right| dx\leq M\int_{\left\{ \left|
u_{n}\right| <t_{0}\right\} }K\left( \left| x\right| \right) \min \left\{
\left| u_{n}\right| ^{q_{1}},\left| u_{n}\right| ^{q_{2}}\right\} dx\leq
M\min \left\{ t_{0}^{q_{1}},t_{0}^{q_{2}}\right\} \left\| K\right\| _{L^{1}(%
\mathbb{R}^{N})}, 
\]
so that, by (\ref{G_pq}), we obtain 
\begin{eqnarray*}
\frac{1}{p}\left\| u_{n}\right\| ^{p}+O\left( 1\right) &=&\int_{\mathbb{R}%
^{N}}G\left( \left| x\right| ,u_{n}\right) dx=\int_{\left\{ \left|
u_{n}\right| <t_{0}\right\} }G\left( \left| x\right| ,u_{n}\right)
dx+\int_{\left\{ \left| u_{n}\right| \geq t_{0}\right\} }G\left( \left|
x\right| ,u_{n}\right) dx \\
&\leq &\tilde{M}\min \left\{ t_{0}^{q_{1}},t_{0}^{q_{2}}\right\} \left\|
K\right\| _{L^{1}(\mathbb{R}^{N})}+\frac{1}{\theta }\int_{\mathbb{R}^{N}}g\left(
\left| x\right| ,u_{n}\right) u_{n}dx+\frac{M}{\theta }\min \left\{
t_{0}^{q_{1}},t_{0}^{q_{2}}\right\} \left\| K\right\| _{L^{1}(\mathbb{R}^{N})}
\\
&=&\left( \tilde{M}+\frac{M}{\theta }\right) \min \left\{
t_{0}^{q_{1}},t_{0}^{q_{2}}\right\} \left\| K\right\| _{L^{1}(\mathbb{R}^{N})}+%
\frac{1}{\theta }\left\| u_{n}\right\| ^{p}+o\left( 1\right) \left\|
u_{n}\right\| .
\end{eqnarray*}
This yields again that $\left\{ \left\| u_{n}\right\| \right\} $ is bounded.
Now, thanks to assumption $\left( \mathcal{S}_{q_{1},q_{2}}^{\prime \prime
}\right) $, we apply Theorem \ref{THM(cpt)} to deduce the existence of $u\in
W_{r}$ such that (up to a subsequence) $u_{n}\rightharpoonup u$ in $W_{r}$
and $u_{n}\rightarrow u$ in $L_{K}^{q_{1}}+L_{K}^{q_{2}}$. Setting 
\[
I_{1}\left( u\right) :=\frac{1}{p}\left\| u\right\| ^{p}\quad \text{and}%
\quad I_{2}\left( u\right) :=I_{1}\left( u\right) -I\left( u\right) 
\]
for brevity, we have that $I_{2}$ is of class $C^{1}$ on $%
L_{K}^{q_{1}}+L_{K}^{q_{2}}$ by \cite[Proposition 3.8]{BPR} and therefore we
get $\left\| u_{n}\right\| ^{p}=I^{\prime }\left( u_{n}\right)
u_{n}+I_{2}^{\prime }\left( u_{n}\right) u_{n}=I_{2}^{\prime }\left(
u\right) u+o\left( 1\right) $. Hence $\lim_{n\rightarrow \infty }\left\|
u_{n}\right\| $ exists and one has $\left\| u\right\| ^{p}\leq
\lim_{n\rightarrow \infty }\left\| u_{n}\right\| ^{p}$ by weak lower
semicontinuity. Moreover, the convexity of $I_{1}:W_{r}\rightarrow \mathbb{R}$
implies 
\[
I_{1}\left( u\right) -I_{1}\left( u_{n}\right) \geq I_{1}^{\prime }\left(
u_{n}\right) \left( u-u_{n}\right) =I^{\prime }\left( u_{n}\right) \left(
u-u_{n}\right) +I_{2}^{\prime }\left( u_{n}\right) \left( u-u_{n}\right)
=o\left( 1\right) 
\]
and thus 
\[
\frac{1}{p}\left\| u\right\| ^{p}=I_{1}\left( u\right) \geq
\lim_{n\rightarrow \infty }I_{1}\left( u_{n}\right) =\frac{1}{p}%
\lim_{n\rightarrow \infty }\left\| u_{n}\right\| ^{p}. 
\]
So $\left\| u_{n}\right\| \rightarrow \left\| u\right\| $ and one concludes
that $u_{n}\rightarrow u$ in $W_{r}$ by the uniform convexity of the norm.%
\endproof%
\bigskip

\noindent \textbf{Proof of Theorem \ref{THM:ex}.}\quad We want to apply the
Mountain-Pass Theorem. To this end, from (\ref{LEM:MPgeom: th}) of Lemma \ref
{LEM:MPgeom} we deduce that, since $L_{0}=0$ and $q_{1},q_{2}>p$, there
exists $\rho >0$ such that 
\begin{equation}
\inf_{u\in W_{r},\,\left\| u\right\| =\rho }I\left( u\right) >0=I\left(
0\right) .  \label{mp-geom}
\end{equation}
Therefore, taking into account Lemma \ref{LEM:PS}, we need only to check
that $\exists \bar{u}\in W_{r}$ such that $\left\| \bar{u}\right\| >\rho $
and $I\left( \bar{u}\right) <0$. To this end, from assumption $\left(
g_{3}\right) $ (which holds in any case, since $\left( g_{1}\right) $ and $%
\left( g_{2}\right) $ imply $\left( g_{3}\right) $), we infer that 
\[
G\left( r,t\right) \geq \frac{G\left( r,t_{0}\right) }{t_{0}^{\theta }}%
t^{\theta }\text{\quad for almost every }r>0\text{~and all }t\geq t_{0}. 
\]
Then, by assumption $\left( \mathbf{V}\right) $, we fix a nonnegative
function $u_{0}\in C_{c}^{\infty }(B_{r_{2}}\setminus \overline{B}%
_{r_{1}})\cap W_{r}$ such that the set $\{x\in \mathbb{R}^{N}:u_{0}\left(
x\right) \geq t_{0}\}$ has positive Lebesgue measure. We now distinguish the
case of assumptions $\left( g_{1}\right) $ and $\left( g_{2}\right) $ from
the case of $K\left( \left| \cdot \right| \right) \in L^{1}(\mathbb{R}^{N})$.
In the first one, $\left( g_{1}\right) $ and $\left( g_{2}\right) $ ensure
that $G\geq 0$ and $G\left( \cdot ,t_{0}\right) >0$ almost everywhere, so
that for every $\lambda >1$ we get 
\begin{eqnarray*}
\int_{\mathbb{R}^{N}}G\left( \left| x\right| ,\lambda u_{0}\right) dx &\geq
&\int_{\left\{ \lambda u_{0}\geq t_{0}\right\} }G\left( \left| x\right|
,\lambda u_{0}\right) dx\geq \frac{\lambda ^{\theta }}{t_{0}^{\theta }}%
\int_{\left\{ \lambda u_{0}\geq t_{0}\right\} }G\left( \left| x\right|
,t_{0}\right) u_{0}^{\theta }dx \\
&\geq &\frac{\lambda ^{\theta }}{t_{0}^{\theta }}\int_{\left\{ u_{0}\geq
t_{0}\right\} }G\left( \left| x\right| ,t_{0}\right) u_{0}^{\theta }dx\geq
\lambda ^{\theta }\int_{\left\{ u_{0}\geq t_{0}\right\} }G\left( \left|
x\right| ,t_{0}\right) dx>0.
\end{eqnarray*}
Since $\theta >p$, this gives 
\[
\lim_{\lambda \rightarrow +\infty }I\left( \lambda u_{0}\right) \leq
\lim_{\lambda \rightarrow +\infty }\left( \frac{\lambda ^{p}}{p}\left\|
u_{0}\right\| ^{p}-\lambda ^{\theta }\int_{\left\{ u_{0}\geq t_{0}\right\}
}G\left( \left| x\right| ,t_{0}\right) dx\right) =-\infty . 
\]
If $K\left( \left| \cdot \right| \right) \in L^{1}(\mathbb{R}^{N})$, assumption 
$\left( g_{3}\right) $ still gives $G\left( \cdot ,t_{0}\right) >0$ almost
everywhere and from (\ref{G_pq}) we infer that 
\[
G\left( r,t\right) \geq -\tilde{M}K\left( r\right) \min \left\{
t_{0}^{q_{1}},t_{0}^{q_{2}}\right\} \text{\quad for almost every }r>0\text{%
~and all }0\leq t\leq t_{0}. 
\]
Therefore, arguing as before about the integral over $\left\{ \lambda
u_{0}\geq t_{0}\right\} $, for every $\lambda >1$ we obtain 
\begin{eqnarray*}
\int_{\mathbb{R}^{N}}G\left( \left| x\right| ,\lambda u_{0}\right) dx
&=&\int_{\left\{ \lambda u_{0}<t_{0}\right\} }G\left( \left| x\right|
,\lambda u_{0}\right) dx+\int_{\left\{ \lambda u_{0}\geq t_{0}\right\}
}G\left( \left| x\right| ,\lambda u_{0}\right) dx \\
&\geq &-\tilde{M}\min \left\{ t_{0}^{q_{1}},t_{0}^{q_{2}}\right\}
\int_{\left\{ \lambda u_{0}<t_{0}\right\} }K\left( \left| x\right| \right)
dx+\lambda ^{\theta }\int_{\left\{ u_{0}\geq t_{0}\right\} }G\left( \left|
x\right| ,t_{0}\right) dx,
\end{eqnarray*}
which implies 
\[
\lim_{\lambda \rightarrow +\infty }I\left( \lambda u_{0}\right) \leq
\lim_{\lambda \rightarrow +\infty }\left( \frac{\lambda ^{p}}{p}\left\|
u_{0}\right\| ^{p}+\tilde{M}\min \left\{ t_{0}^{q_{1}},t_{0}^{q_{2}}\right\}
\left\| K\right\| _{L^{1}(\mathbb{R}^{N})}-\lambda ^{\theta }\int_{\left\{
u_{0}\geq t_{0}\right\} }G\left( \left| x\right| ,t_{0}\right) dx\right)
=-\infty . 
\]
So, in any case, we can take $\bar{u}=\lambda u_{0}$ with $\lambda $
sufficiently large and the Mountain-Pass Theorem provides the existence of a
nonzero critical point $u\in W_{r}$ for $I$. Since (\ref{g=0}) implies $%
I^{\prime }\left( u\right) u_{-}=-\left\| u_{-}\right\| ^{p}$ (where $%
u_{-}\in W_{r}$ is the negative part of $u$), one concludes that $u_{-}=0$,
i.e., $u$ is nonnegative.%
\endproof%
\bigskip

In order to conclude the proof of Theorem \ref{THM:mult}, we recall from 
\cite[Corollary 2.19]{BPR} that for every $p_{1},p_{2}\in \left( 1,+\infty
\right) $ and $u\in L_{K}^{p_{1}}+L_{K}^{p_{2}}$ one has 
\begin{equation}
\left\| u\right\| _{L_{K}^{p_{1}}+L_{K}^{p_{2}}}\leq \left\| u\right\|
_{L_{K}^{\min \left\{ p_{1},p_{2}\right\} }(\Lambda _{u})}+\left\| u\right\|
_{L_{K}^{\max \left\{ p_{1},p_{2}\right\} }(\Lambda _{u}^{c})},\quad \text{%
where}\quad \Lambda _{u}:=\left\{ x\in \mathbb{R}^{N}:\left| u\left( x\right)
\right| >1\right\} .  \label{Vu}
\end{equation}

\noindent \textbf{Proof of Theorem \ref{THM:mult}.}\quad By the oddness
assumption $\left( g_{5}\right) $, one has $I\left( u\right) =I\left(
-u\right) $ for all $u\in W_{r}$ and thus we can apply the Symmetric
Mountain-Pass Theorem (see e.g. \cite[Chapter 1]{Rabi}). To this end, we
deduce (\ref{mp-geom}) as in the proof of Theorem \ref{THM:ex} and
therefore, thanks to Lemma \ref{LEM:PS}, we need only to show that $I$
satisfies the following geometrical condition: for any finite dimensional
subspace $X\neq \left\{ 0\right\} $ of $W_{r}$ there exists $R>0$ such that $%
I\left( u\right) \leq 0$ for all $u\in X$ with $\left\| u\right\| \geq R$.
In fact, it is sufficient to prove that any diverging sequence in $X$ admits
a subsequence on which $I$ is nonpositive. So, let $\left\{ u_{n}\right\}
\subseteq X$ be such that $\left\| u_{n}\right\| \rightarrow +\infty $.
Since all norms are equivalent on $X$, by (\ref{Vu}) one has 
\begin{equation}
\left\| u_{n}\right\| _{L_{K}^{p_{1}}(\Lambda _{u_{n}})}+\left\|
u_{n}\right\| _{L_{K}^{p_{2}}(\Lambda _{u_{n}}^{c})}\geq \left\|
u_{n}\right\| _{L_{K}^{q_{1}}+L_{K}^{q_{2}}}\geq m_{1}\left\| u_{n}\right\|
\rightarrow +\infty  \label{THM2-pf: eq_norms}
\end{equation}
for some constant $m_{1}>0$, where $p_{1}:=\min \left\{ q_{1},q_{2}\right\} $
and $p_{2}:=\max \left\{ q_{1},q_{2}\right\} $. Hence, up to a subsequence,
at least one of the sequences $\{\left\| u_{n}\right\|
_{L_{K}^{p_{1}}(\Lambda _{u_{n}})}\}$, $\{\left\| u_{n}\right\|
_{L_{K}^{p_{2}}(\Lambda _{u_{n}}^{c})}\}$ diverges. We now use assumptions $%
\left( g_{4}\right) $ and $\left( g_{5}\right) $ to deduce that 
\[
G\left( r,t\right) \geq mK\left( r\right) \min \left\{ \left| t\right|
^{q_{1}},\left| t\right| ^{q_{2}}\right\} \quad \text{for almost every }r>0~%
\text{and all }t\in \mathbb{R}, 
\]
which implies 
\[
\int_{\mathbb{R}^{N}}G\left( \left| x\right| ,u_{n}\right) dx\geq
m\int_{\Lambda _{u_{n}}}K\left( \left| x\right| \right) \left| u_{n}\right|
^{p_{1}}dx+m\int_{\Lambda _{u_{n}}^{c}}K\left( \left| x\right| \right)
\left| u_{n}\right| ^{p_{2}}dx. 
\]
Hence, using inequalities (\ref{THM2-pf: eq_norms}), there exists a constant 
$m_{2}>0$ such that 
\[
I\left( u_{n}\right) \leq m_{2}\left( \left\| u_{n}\right\|
_{L_{K}^{p_{1}}(\Lambda _{u_{n}})}^{p}+\left\| u_{n}\right\|
_{L_{K}^{p_{2}}(\Lambda _{u_{n}}^{c})}^{p}\right) -m\left( \left\|
u_{n}\right\| _{L_{K}^{p_{1}}(\Lambda _{u_{n}})}^{p_{1}}+\left\|
u_{n}\right\| _{L_{K}^{p_{2}}(\Lambda _{u_{n}}^{c})}^{p_{2}}\right) , 
\]
so that $I\left( u_{n}\right) \rightarrow -\infty $ since $p_{1},p_{2}>p$.
The Symmetric Mountain-Pass Theorem thus implies the existence of an
unbounded sequence of critical values for $I$ and this completes the proof.%
\endproof%
\bigskip

\begin{lem}
\label{LEM:bdd}Under the assumptions of each of Theorems \ref{THM:ex-sub}
and \ref{THM:mult-sub}, the functional $I:W_{r}\rightarrow \mathbb{R}$ is
bounded from below and coercive. In particular, if $g$ satisfies $\left(
g_{6}\right) $, then 
\begin{equation}
\inf_{v\in W_{r}}I\left( v\right) <0.  \label{LEM:bdd: inf<0}
\end{equation}
\end{lem}

\proof%
The fact that $I$ is bounded below and coercive on $W_{r}$ is a consequence
of Lemma \ref{LEM:MPgeom}. Indeed, the result readily follows from (\ref
{LEM:MPgeom: th}) if $q_{1},q_{2}\in \left( 1,p\right) $, while, if $\max
\left\{ q_{1},q_{2}\right\} =p>\min \left\{ q_{1},q_{2}\right\} >1$, we fix $%
\varepsilon <1/p$ and use the second part of the lemma in order to get 
\[
I\left( u\right) \geq \left( \frac{1}{p}-\varepsilon \right) \left\|
u\right\| ^{p}-c\left( \varepsilon \right) \left\| u\right\| ^{\min \left\{
q_{1},q_{2}\right\} }-L_{0}\left\| u\right\| \qquad \text{for all }u\in
W_{r}, 
\]
which yields again the conclusion. In order to prove (\ref{LEM:bdd: inf<0})
under assumption $\left( g_{6}\right) $, we use assumption $\left( \mathbf{V}%
\right) $ to fix a function $u_{0}\in C_{c}^{\infty }(B_{r_{2}}\setminus 
\overline{B}_{r_{1}})\cap W_{r}$ such that $0\leq u_{0}\leq t_{0}$, $%
u_{0}\neq 0$. Then, by $\left( g_{6}\right) $, for every $0<\lambda <1$ we
get that $\lambda u_{0}\in W_{r}$ satisfies 
\[
I\left( \lambda u_{0}\right) =\frac{1}{p}\left\| \lambda u_{0}\right\|
^{p}-\int_{\mathbb{R}^{N}}G\left( \left| x\right| ,\lambda u_{0}\right) dx\leq 
\frac{\lambda ^{2}}{2}\left\| u_{0}\right\| ^{2}-\lambda ^{\theta }m\int_{%
\mathbb{R}^{N}}K\left( \left| x\right| \right) u_{0}^{\theta }dx. 
\]
Since $\theta <p$, this implies $I\left( \lambda u_{0}\right) <0$ for $%
\lambda $ sufficiently small and therefore (\ref{LEM:bdd: inf<0}) ensues.%
\endproof%
\bigskip

\noindent \textbf{Proof of Theorem \ref{THM:ex-sub}.}\quad Let 
\[
\mu :=\inf_{v\in W_{r}}I\left( v\right) 
\]
and take any minimizing sequence $\left\{ v_{n}\right\} $ for $\mu $. From
Lemma \ref{LEM:bdd} we have that the functional $I:W_{r}\rightarrow \mathbb{R}$
is bounded from below and coercive, so that $\mu \in \mathbb{R}$ and $\left\{
v_{n}\right\} $ is bounded in $W_{r}$. Thanks to Theorem \ref{THM(cpt)} and
assumption $\left( \mathcal{S}_{q_{1},q_{2}}^{\prime \prime }\right) $, the
embedding $W_{r}\hookrightarrow L_{K}^{q_{1}}+L_{K}^{q_{2}}$ is compact and
thus we can assume that there exists $u\in W_{r}$ such that, up to a
subsequence, one has $v_{n}\rightharpoonup u$ in $W_{r}$ and $%
v_{n}\rightarrow u$ in $L_{K}^{q_{1}}+L_{K}^{q_{2}}$. Then, thanks to $%
\left( g_{0}\right) $ and the continuity of the functional $v\mapsto \int_{%
\mathbb{R}^{N}}\left( G\left( \left| x\right| ,v\right) -g\left( \left|
x\right| ,0\right) v\right) dx$ on $L_{K}^{q_{1}}+L_{K}^{q_{2}}$ (which
follows from $\left( g\right) $, $\left( f_{q_{1},q_{2}}\right) $ and 
\cite[Proposition 3.8]{BPR}), $u$ satisfies 
\[
\int_{\mathbb{R}^{N}}G\left( \left| x\right| ,v_{n}\right) dx=\int_{\mathbb{R}%
^{N}}\left( G\left( \left| x\right| ,v_{n}\right) -g\left( \left| x\right|
,0\right) v_{n}\right) dx+\int_{\mathbb{R}^{N}}g\left( \left| x\right|
,0\right) v_{n}dx\rightarrow \int_{\mathbb{R}^{N}}G\left( \left| x\right|
,u\right) dx. 
\]
By the weak lower semi-continuity of the norm, this implies 
\[
I\left( u\right) =\frac{1}{p}\left\| u\right\| ^{p}-\int_{\mathbb{R}%
^{N}}G\left( \left| x\right| ,u\right) dx\leq \lim_{n\rightarrow \infty
}\left( \frac{1}{p}\left\| v_{n}\right\| ^{p}-\int_{\mathbb{R}^{N}}G\left(
\left| x\right| ,v_{n}\right) dx\right) =\mu 
\]
and thus we conclude $I\left( u\right) =\mu $. It remains to show that $%
u\neq 0$. This is obvious if $g$ satisfies $\left( g_{6}\right) $, since $%
\mu <0$ by Lemma \ref{LEM:bdd}. If $\left( g_{7}\right) $ holds, assume by
contradiction that $u=0$. Since $u$ is a critical point of $I\in C^{1}(W_{r};%
\mathbb{R})$, from (\ref{PROP:diff: I'(u)h=}) we get 
\[
\int_{\mathbb{R}^{N}}g\left( \left| x\right| ,0\right) h\,dx=0,\quad \forall
h\in C_{c,\mathrm{rad}}^{\infty }(B_{r_{2}}\setminus \overline{B}%
_{r_{1}})\subset W_{r}. 
\]
This implies $g\left( \cdot ,0\right) =0$ almost everywhere in $\left(
r_{1},r_{2}\right) $, which is a contradiction.%
\endproof%
\bigskip

\noindent \textbf{Proof of Corollary \ref{COR:ex-sub}.}\quad Setting 
\[
\widetilde{g}\left( r,t\right) :=\left\{ 
\begin{array}{lll}
g\left( r,t\right) &  & \text{if }t\geq 0 \\ 
2g\left( r,0\right) -g\left( r,\left| t\right| \right) &  & \text{if }t<0
\end{array}
\right. 
\]
and 
\[
\widetilde{G}\left( r,t\right) :=\int_{0}^{t}\widetilde{g}\left( r,s\right)
ds=\left\{ 
\begin{array}{lll}
G\left( r,t\right) &  & \text{if }t\geq 0 \\ 
2g\left( r,0\right) t+G\left( r,\left| t\right| \right) &  & \text{if }t<0,
\end{array}
\right. 
\]
it is easy to check that the function $\widetilde{g}$ still satisfies all
the assumptions of Theorem \ref{THM:ex-sub}. Then there exists $\widetilde{u}%
\neq 0$ such that 
\[
\widetilde{I}\left( \widetilde{u}\right) =\min_{u\in W_{r}}\widetilde{I}%
\left( u\right) ,\quad \text{where}\quad \widetilde{I}\left( u\right) :=%
\frac{1}{p}\left\| u\right\| ^{p}-\int_{\mathbb{R}^{N}}\widetilde{G}\left(
\left| x\right| ,u\right) dx. 
\]
For every $u\in W_{r}$ one has 
\begin{eqnarray}
\widetilde{I}\left( u\right) &=&\frac{1}{p}\left\| u\right\|
^{p}-\int_{\left\{ u\geq 0\right\} }G\left( \left| x\right| ,u\right)
dx-2\int_{\left\{ u<0\right\} }g\left( \left| x\right| ,0\right)
u\,dx-\int_{\left\{ u<0\right\} }G\left( \left| x\right| ,\left| u\right|
\right) dx  \nonumber \\
&=&\frac{1}{p}\left\| u\right\| ^{p}-\int_{\mathbb{R}^{N}}G\left( \left|
x\right| ,\left| u\right| \right) dx+2\int_{\mathbb{R}^{N}}g\left( \left|
x\right| ,0\right) u_{-}\,dx  \label{I-tilde=} \\
&=&I\left( \left| u\right| \right) +2\int_{\mathbb{R}^{N}}g\left( \left|
x\right| ,0\right) u_{-}\,dx,  \nonumber
\end{eqnarray}
which implies that $\widetilde{u}$ satisfies (\ref{COR:ex-sub: th}), as one
readily checks that 
\[
\inf_{u\in W_{r}}\left( I\left( \left| u\right| \right) +2\int_{\mathbb{R}%
^{N}}g\left( \left| x\right| ,0\right) u_{-}\,dx\right) =\inf_{u\in
W_{r},\,u\geq 0}I\left( u\right) . 
\]
Moreover, since $G\left( r,\left| t\right| \right) =\widetilde{G}\left(
r,\left| t\right| \right) $ and $g\left( \cdot ,0\right) \geq 0$, (\ref
{I-tilde=}) gives $\widetilde{I}\left( u\right) \geq \widetilde{I}\left(
\left| u\right| \right) $ for every $u\in W_{r}$ and hence $\left| 
\widetilde{u}\right| \in W_{r}$ is still a minimizer for $\widetilde{I}$, so
that we can assume $\widetilde{u}\geq 0$. Finally, $\widetilde{u}$ is a
critical point for $I$ since $\widetilde{u}$ is a critical point of $%
\widetilde{I}$ and $\widetilde{g}\left( r,t\right) =g\left( r,t\right) $ for
every $t\geq 0$. 
\endproof%
\bigskip

In proving Theorem \ref{THM:mult-sub} we will use a well known abstract
result, which we recall it here in a version from \cite{Wang01}.

\begin{thm}[{\cite[Lemma 2.4]{Wang01}}]
\label{THM:wang}Let $X$ be a real Banach space and let $J\in C^{1}(X;\mathbb{R}%
) $. Assume that $J$ satisfies the Palais-Smale condition, is even, bounded
from below and such that $J\left( 0\right) =0$. Assume furthermore that $%
\forall k\in \mathbb{N}\setminus \left\{ 0\right\} $ there exist $\rho _{k}>0$
and a $k$-dimensional subspace $X_{k}$ of $X$ such that 
\begin{equation}
\sup_{u\in X_{k},\,\left\| u\right\| _{X}=\rho _{k}}J\left( u\right) <0.
\label{THM:wang: geom}
\end{equation}
Then $J$ has a sequence of critical values $c_{k}<0$ such that $%
\displaystyle%
\lim_{k\rightarrow \infty }c_{k}=0$.
\end{thm}

\noindent 

\noindent \textbf{Proof of Theorem \ref{THM:mult-sub}.}\quad Since $%
I:W_{r}\rightarrow \mathbb{R}$ is even by assumption $\left( g_{5}\right) $ and
bounded below by Lemma \ref{LEM:bdd}, for applying Theorem \ref{THM:wang}
(with $X=W_{r}$ and $J=I$) we need only to show that $I$ satisfies the
Palais-Smale condition and the geometric condition (\ref{THM:wang: geom}).
By coercivity (see Lemma \ref{LEM:bdd} again), every Palais-Smale sequence
for $I$ is bounded in $W_{r}$ and one obtains the existence of a strongly
convergent subsequence as in the proof of Lemma \ref{LEM:PS}. In order to
check (\ref{THM:wang: geom}), we first deduce from $\left( g_{5}\right) $
and $\left( g_{6}\right) $ that 
\begin{equation}
G\left( r,t\right) \geq mK\left( r\right) \left| t\right| ^{\theta }\quad 
\text{for almost every }r>0~\text{and all }\left| t\right| \leq t_{0}.
\label{G>}
\end{equation}
Then, for any $k\in \mathbb{N}\setminus \left\{ 0\right\} $, we take $k$
linearly independent functions $\phi _{1},...,\phi _{k}\in C_{c,\mathrm{rad}%
}^{\infty }(B_{r_{2}}\setminus \overline{B}_{r_{1}})$ such that $0\leq \phi
_{i}\leq t_{0}$ for every $i=1,...,k$ and set 
\[
X_{k}:=\linspan\left\{ \phi _{1},...,\phi _{k}\right\} \text{\quad
and\quad }\left\| \lambda _{1}\phi _{1}+...+\lambda _{k}\phi _{k}\right\|
_{X_{k}}:=\max_{1\leq i\leq k}\left| \lambda _{i}\right| . 
\]
This defines a subspace of $W_{r}$ by assumption $\left( \mathbf{V}\right) $
and all the norms are equivalent on $X_{k}$, so that there exist $%
m_{k},l_{k}>0$ such that for all $u\in X_{k}$ one has 
\begin{equation}
\left\| u\right\| _{X_{k}}\leq m_{k}\left\| u\right\| \quad \text{and}\quad
\left\| u\right\| _{L_{K}^{\theta }(\mathbb{R}^{N})}^{\theta }\geq l_{k}\left\|
u\right\| ^{\theta }.  \label{equivalent}
\end{equation}
Fix $\rho _{k}>0$ small enough that $km_{k}\rho _{k}<1$ and $\rho
_{k}^{p}/p-ml_{k}\rho _{k}^{\theta }<0$ (which is possible since $\theta <p$%
) and take any $u=\lambda _{1}\phi _{1}+...+\lambda _{k}\phi _{k}\in X_{k}$
such that $\left\| u\right\| =\rho _{k}$. Then by (\ref{equivalent}) we have 
$\left| \lambda _{i}\right| \leq \left\| u\right\| _{X_{k}}\leq m_{k}\rho
_{k}<1/k$ for every $i=1,...,k$ and therefore 
\[
\left| u\left( x\right) \right| \leq \sum_{i=1}^{k}\left| \lambda
_{i}\right| \phi _{i}\left( x\right) \leq t_{0}\sum_{i=1}^{k}\left| \lambda
_{i}\right| <t_{0}\qquad \text{for all }x\in \mathbb{R}^{N}. 
\]
By (\ref{G>}) and (\ref{equivalent}), this implies 
\[
\int_{\mathbb{R}^{N}}G\left( \left| x\right| ,u\right) dx\geq m\int_{\mathbb{R}%
^{N}}K\left( \left| x\right| \right) \left| u\right| ^{\theta }dx\geq
ml_{k}\left\| u\right\| ^{\theta } 
\]
and hence we get $I\left( u\right) \leq \left\| u\right\|
^{p}/p-ml_{k}\left\| u\right\| ^{\theta }=\rho _{k}^{p}/p-ml_{k}\rho
_{k}^{\theta }<0$. This proves (\ref{THM:wang: geom}) and the conclusion
thus follows from Theorem \ref{THM:wang}.%
\endproof%


\end{document}